\documentclass[11pt]{article}
\usepackage[mathscr]{eucal}
\usepackage{mathrsfs}
\usepackage{amssymb}
\usepackage{calligra}
\usepackage{fontenc}
\usepackage{amsbsy} 
\usepackage{graphicx}
\graphicspath{%
    {converted_graphics/}
    {/}
}
\begin{document}

\newcommand{\simge}{\ba{cc}\vspace*{-2.4mm}>\\ \sim\ea }
\newcommand{\simle}{\ba{cc}\vspace*{-2.4mm}<\\ \sim\ea }
\newcommand{\Cdot}{\!\cdot\!}
\newcommand{\sq}{{$\sqcap\!\!\!\!\sqcup$}}
\newcommand{\Eu}{{\rm I\,\!\! E}}
\newcommand{\Io}{\Int{\Omega}{}}
\newcommand{\Id}{\Int{\cald}{}}
\newcommand{\Div}{\mbox{\rm div}\,}
\newcommand{\tr}{\mbox{\rm tr}\,}
\newcommand{\grad}{\mbox{\rm grad}\,}
\newcommand{\supp}{\mbox{\rm supp}\,}
\newcommand{\curl}{\mbox{\rm curl}\,}
\newcommand{\Ido}{\Int{\partial\Omega}{}}
\newcommand{\IdS}{\Int{\Sigma}{}}
\newcommand{\Oint}[2]{{\displaystyle \oint_{#1}^{#2}}}
\newcommand{\Int}[2]{{\displaystyle \int_{ #1}^{ #2}}}
\newcommand{\Lim}[1]{{\displaystyle \lim_{ #1}}}
\newcommand{\Limsup}[1]{{\displaystyle \limsup_{\footnotesize #1}}}
\newcommand{\Liminf}[1]{{\displaystyle \liminf_{\footnotesize #1}}}
\newcommand{\Sup}[1]{{\displaystyle \sup_{#1}}}
\newcommand{\Inf}[1]{{\displaystyle \inf_{#1}}}
\newcommand{\Max}[1]{{\displaystyle \max_{#1}}}
\newcommand{\Min}[1]{{\displaystyle \min_{#1}}}
\newcommand{\Sum}[2]{{\displaystyle \sum_{#1}^{#2}}}
\newcommand{\Prod}[2]{{\displaystyle \prod_{#1}^{#2}}}
\newcommand{\BCup}[2]{{\displaystyle \bigcup_{#1}^{#2}}}
\newcommand{\BCap}[2]{{\displaystyle \bigcap_{#1}^{#2}}}
\newcommand{\Frac}[2]{\displaystyle{\frac{\displaystyle{#1}}{\displaystyle{#2}}}}
\newcommand{\norm}[1]{\left\|{#1}\right\|}
\newcommand{\Norm}[1]{\langle\langle{#1}\rangle\rangle_q}
\newcommand{\No}[1]{\langle\!\langle{#1}\rangle\!\rangle}
\newcommand{\NO}[1]{{\langle{#1}\rangle}_{\lambda,q}}
\newcommand{\beea}{\begin{eqnarray}}
\newcommand{\eeea}{\end{eqnarray}}
\newcommand{\ms}{\medskip\smallskip}
\newcommand{\bs}{\bigskip}
\newcommand{\ps}{\par\smallskip}
\newcommand{\bfe}{{\mbox{\boldmath $e$}} }
\newcommand{\pni}{{\par\noindent}}
\newcommand{\bfq}{{\mbox{\boldmath $q$}} }
\newcommand{\bfz}{{\mbox{\boldmath $z$}} }
\newcommand{\0}{{\mbox{\boldmath $0$}} }
\newcommand{\LE}{\!\!\!&\le&\!\!\!}
\newcommand{\BL}[1]{{\par\smallskip{\bf Lemma #1.}}}
\newcommand{\BT}[1]{{\par\smallskip{\bf Theorem #1.}}}
\newcommand{\Ln}{[\!|}
\newcommand{\Rn}{|\!]}
\newcommand{\n}[1]{{\Ln{#1}\Rn}} 
\newcommand{\nq}[1]{{\Ln{#1}\Rn}_{q}} 
\newcommand{\nqr}[1]{{\Ln{#1}\Rn}_{q,r}} 
\newcommand{\Nq}[1]{{\langle{#1}\rangle}_{q}} 
\newcommand{\Nql}[1]{{\langle{#1}\rangle}_{\lambda,q}} 
\newcommand{\Nqr}[1]{{\langle{#1}\rangle}_{q,r}}
\newcommand{\N}[1]{{|\!\!|\!\!|\,{#1}\,|\!|\!\!|_2}}
\newcommand{\EA}[2]{$$#1$$%
\vspace{-6.mm}
\begin{equation}
\end{equation}
\vspace{-6.mm}
$$
#2
\setlength{\belowdisplayskip}{3mm}
\setlength{\belowdisplayshortskip}{3mm}
$$
}
\newcommand{\A}[2]{$$#1$$%
\vspace{-4.mm}
$$
#2
\setlength{\belowdisplayskip}{3mm}
\setlength{\belowdisplayshortskip}{3mm}
$$
}
\newcommand{\BF}{\begin{footnotesize}}
\newcommand{\EF}{\end{footnotesize}}
\setlength{\jot}{.15in}
\newcommand{\pde}[2]{{\displaystyle \frac{\mbox{$\partial #1$}}{\mbox{$\partial #2$}}}}
\newcommand{\ode}[2]{{\displaystyle \frac{\mbox{$d #1$}}{\mbox{$d #2$}}}}
\newcommand{\f}[2]{\frac{\mbox{$#1$}}{\mbox{$ #2$}}}
\newcommand{\bi}{\begin{itemize}}
\newcommand{\ei}{\end{itemize}}
\newcommand{\ed}{\end{document}}
\newcommand{\be}{\begin{equation}}
\newcommand{\ba}{\begin{array}}
\newcommand{\ea}{\end{array}}
\newcommand{\ee}{\end{equation}}
\newcommand{\eeq}[1]{\label{eq:#1}\end{equation}}
\newcommand{\real}{{\mathbb R}}
\newcommand{\compl}{{\mathbb C}}
\def\Id{\mbox{\boldmath $1$}}
\def\zero{\mbox{\boldmath $0$}}
\newcommand{\PP}{{\rm I\!\!\,P}}
\newcommand{\nat}{{\mathbb N}}
\newcommand{\bfpsi}{\mbox{\boldmath $\psi$}}
\newcommand{\bfchi}{\mbox{\boldmath $\chi$}}
\newcommand{\bfomega}{\mbox{\boldmath $\omega$}}
\newcommand{\bfvaromega}{\mbox{\boldmath $\varpi$}}
\newcommand{\bfOmega}{\mbox{\boldmath $\Omega$}}
\newcommand{\bfTheta}{\mbox{\boldmath $\Theta$}}
\newcommand{\bfxi}{\mbox{\boldmath $\xi$}}
\newcommand{\bfmu}{\mbox{\boldmath $\mu$}}
\newcommand{\bfx}{\mbox{\boldmath $x$}}
\newcommand{\bfy}{\mbox{\boldmath $y$}}
\newcommand{\bfPsi}{\mbox{\boldmath $\Psi$}}
\newcommand{\bfphi}{\mbox{\boldmath $\varphi$}}
\newcommand{\bfhi}{\mbox{\boldmath $\phi$}}
\newcommand{\bfPhi}{\mbox{\boldmath $\Phi$}}
\newcommand{\bfv}{{\mbox{\boldmath $v$}} }
\newcommand{\bfu}{{\mbox{\boldmath $u$}} }
\newcommand{\bfsf}{{\mbox{\footnotesize\boldmath $s$}} }
\newcommand{\bfuf}{{\mbox{\footnotesize\boldmath $u$}} }
\newcommand{\bfw}{{\mbox{\boldmath $w$}} }
\newcommand{\bff}{{\mbox{\boldmath $f$}} }
\newcommand{\bfa}{{\mbox{\boldmath $a$}} }
\newcommand{\bfi}{{\mbox{\boldmath $i$}} }
\newcommand{\bfj}{{\mbox{\boldmath $j$}} }
\newcommand{\bfc}{{\mbox{\boldmath $c$}} }
\newcommand{\bfo}{{\mbox{\boldmath $o$}} }
\newcommand{\bfp}{{\mbox{\boldmath $p$}} }
\newcommand{\bfkp}{{\mbox{\footnotesize{\boldmath $k$}}} }
\newcommand{\bfka}{{\mbox{\footnotesize{\boldmath $k^*$}}} }
\newcommand{\bft}{{\mbox{\boldmath $t$}} }
\newcommand{\bfd}{{\mbox{\boldmath $d$}} }
\newcommand{\bfl}{{\mbox{\boldmath $l$}} }
\newcommand{\bfr}{{\mbox{\boldmath $r$}} }
\newcommand{\bfk}{{\mbox{\boldmath $k$}} }
\newcommand{\bfA}{{\mbox{\boldmath $A$}} }
\newcommand{\bfS}{{\mbox{\boldmath $S$}} }
\newcommand{\bfO}{{\mbox{\boldmath $O$}} }
\newcommand{\bfM}{{\mbox{\boldmath $M$}} }
\newcommand{\bfP}{{\mbox{\boldmath $P$}} }
\newcommand{\bfB}{{\mbox{\boldmath $B$}} }
\newcommand{\bfR}{{\mbox{\boldmath $R$}} }
\newcommand{\bfC}{{\mbox{\boldmath $C$}} }
\newcommand{\bfD}{{\mbox{\boldmath $D$}} }
\newcommand{\bfQ}{{\mbox{\boldmath $Q$}} }
\newcommand{\bfZ}{{\mbox{\boldmath $Z$}} }
\newcommand{\bfG}{{\mbox{\boldmath $G$}} }
\newcommand{\bfE}{{\mbox{\boldmath $E$}} }
\newcommand{\bfX}{{\mbox{\boldmath $X$}} }
\newcommand{\bfY}{{\mbox{\boldmath $Y$}} }
\newcommand{\bfH}{{\mbox{\boldmath $H$}} }
\newcommand{\bfI}{{\mbox{\boldmath $I$}} }
\newcommand{\bfJ}{{\mbox{\boldmath $J$}} }
\newcommand{\bfN}{{\mbox{\boldmath $N$}} }
\newcommand{\bfh}{{\mbox{\boldmath $h$}} }
\newcommand{\bfm}{{\mbox{\boldmath $m$}} }
\newcommand{\bfone}{{\mbox{\boldmath $1$}} }
\newcommand{\hs}{{\rm I}\!\!\,{\rm R}^3_+}
\newcommand{\cala}{{\cal A}}
\newcommand{\calb}{{\cal B}}
\newcommand{\calc}{{\cal C}}
\newcommand{\cald}{{\cal D}}
\newcommand{\cale}{{\cal E}}
\newcommand{\calf}{{\cal F}}
\newcommand{\calg}{{\cal G}}
\newcommand{\calh}{{\cal H}}
\newcommand{\cali}{{\cal I}}
\newcommand{\calj}{{\cal J}}
\newcommand{\calk}{{\cal K}}
\newcommand{\call}{{\cal L}}
\newcommand{\calm}{{\cal M}}
\newcommand{\caln}{{\cal N}}
\newcommand{\calo}{{\cal O}}
\newcommand{\calp}{{\cal P}}
\newcommand{\calq}{{\cal Q}}
\newcommand{\calr}{{\cal R}}
\newcommand{\cals}{{\cal S}}
\newcommand{\calt}{{\cal T}}
\newcommand{\calu}{{\cal U}}
\newcommand{\calv}{{\cal V}}
\newcommand{\calx}{{\cal X}}
\newcommand{\caly}{{\cal Y}}
\newcommand{\calw}{{\cal W}}
\newcommand{\calz}{{\cal Z}}
\newcommand{\bfsigma}{\mbox{\boldmath $\sigma$}}
\newcommand{\bfSigma}{\mbox{\boldmath $\Sigma$}}
\newcommand{\bftau}{\mbox{\boldmath $\tau$}}
\newcommand{\bfeta}{\mbox{\boldmath $\eta$}}
\newcommand{\bfT}{{\mbox{\boldmath $T$}} }
\newcommand{\bfV}{{\mbox{\boldmath $V$}} }
\newcommand{\bfU}{{\mbox{\boldmath $U$}} }
\newcommand{\bfW}{{\mbox{\boldmath $W$}} }
\newcommand{\bfF}{{\mbox{\boldmath $F$}} }
\newcommand{\bfK}{{\mbox{\boldmath $K$}} }
\newcommand{\bfL}{{\mbox{\boldmath $L$}} }
\newcommand{\bfb}{{\mbox{\boldmath $b$}} }
\newcommand{\bfg}{{\mbox{\boldmath $g$}} }
\newcommand{\bfn}{{\mbox{\boldmath $n$}} }
\newcommand{\bfs}{{\mbox{\boldmath $s$}} }
\newcommand{\cf}{{\it cf.} }
\newcommand{\io}{\int_\Omega}
\newcommand{\1}{\item[({\it i})]}
\newcommand{\2}{\item[({\it ii})]}
\newcommand{\3}{\item[({\it iii})]}
\newcommand{\4}{\item[({\it iv})]}
\newcommand{\5}{\item[({\it v})]}
\newcommand{\6}{\item[({\it vi})]}
\newcommand{\7}{\item[({\it vii})]}
\newcommand{\8}{\item[({\it viii})]}
\newcommand{\9}{\item[({\it xi})]}
\newcommand{\ido}{\int_{\partial\Omega}}
\newcommand{\half}{\mbox{$\frac{1}{2}$}}
\def\parallel{\|}
\def\mid{|}
\def\Bbb R{\real}
\def\hat{\widehat}
\def\tilde{\widetilde}
\def\bar{\overline}
\newcommand{\threehalves}{3\over 2}
\newcommand{\bfPi}{\mbox{\boldmath $\Pi$}}
\newcommand{\bfXi}{\mbox{\boldmath $\Xi$}}
\newcommand{\bfalpha}{\mbox{\boldmath $\alpha$}}
\newcommand{\bfbeta}{\mbox{\boldmath $\beta$}}
\newcommand{\bfgamma}{\mbox{\boldmath $\gamma$}}
\newcommand{\bfdelta}{\mbox{\boldmath $\delta$}}
\newcommand{\bfzeta}{\mbox{\boldmath $\zeta$}}
\newcommand{\bfUpsilon}{\mbox{\boldmath $\Upsilon$}}
\newcommand{\bfGamma}{\mbox{\boldmath $\Gamma$}}
\newcommand{\bfcala}{\mbox{\boldmath ${\cal A}$}}
\newcommand{\bfcalm}{\mbox{\boldmath ${\cal M}$}}
\newcommand{\bfcaln}{\mbox{\boldmath ${\cal N}$}}
\newcommand{\bfcalq}{\mbox{\boldmath ${\cal Q}$}}
\newcommand{\bfcalb}{\mbox{\boldmath ${\cal B}$}}
\newcommand{\bfcalc}{\mbox{\boldmath ${\cal C}$}}
\newcommand{\bfcali}{\mbox{\boldmath ${\cal I}$}}
\newcommand{\bfcalg}{\mbox{\boldmath ${\cal G}$}}
\newcommand{\bfcalh}{\mbox{\boldmath ${\cal H}$}}
\newcommand{\bfcalk}{\mbox{\boldmath ${\cal K}$}}
\newcommand{\bfcalt}{\mbox{\boldmath ${\cal T}$}}
\newcommand{\bfcalx}{\mbox{\boldmath ${\cal X}$}}
\newcommand{\bfcall}{\mbox{\boldmath ${\cal L}$}}
\newcommand{\bfcalf}{\mbox{\boldmath ${\cal F}$}}
\newcommand{\bfcalr}{\mbox{\boldmath ${\cal R}$}}
\newcommand{\bfcals}{\mbox{\boldmath ${\cal S}$}}
\newcommand{\bfcalw}{\mbox{\boldmath ${\cal W}$}}
\newcommand{\bfcalu}{\mbox{\boldmath ${\cal U}$}}
\newcommand{\bfcalv}{\mbox{\boldmath ${\cal V}$}}
\newcommand{\bfcalz}{\mbox{\boldmath ${\cal Z}$}}
\pagenumbering{roman}
\newcommand{\art}[6]{{\I[{\sc #1,}] {#2}, {\it #3}, {\bf #4}, {#5} {[#6]}}}
\newcommand{\ED}{\end{description}}
\newcommand{\I}{\item }
\newcommand{\ra}{\rm a}
\newcommand{\rb}{\rm b}
\newcommand{\rc}{\rm c}
\newcommand{\Hsp}{{\rm I}\!\!\,{\rm R}^n_+}
\newcommand{\Hsn}{{\rm I}\!\!\,{\rm R}^n_-}
\newcommand{\po}[1]{\mbox{$\displaystyle \frac{\mbox{$\partial #1$}}
{\mbox{$\partial x_{1}$}}$}}
\newcommand{\PO}[1]{\mbox{$\displaystyle \frac{\mbox{$\partial #1$}}
{\mbox{$\partial y_{1}$}}$}}
\newcommand{\OP}{\left(\Delta+2\lambda\PO{}\right)}
\newcommand{\op}{\left(\Delta+2\lambda\po{}\right)}
\newcommand{\ft}[1]{
\Frac{1}{(2\pi)^{n/2}}\Int{{\Bbb R}^{n}}{}e^{i{\bf x}\cdot \bfxi}
#1(\xi)d\xi}
\newcommand{\Ft}[1]{
\Frac{1}{2\pi}\Int{{\Bbb R}^{2}}{}e^{i{x}\cdot \xi}
#1(\xi)d\xi}
\newcommand{\Z}{\item[({\it a})]}
\newcommand{\B}{\item[({\it b})]}
\newcommand{\C}{\item[({\it c})]}
\newcommand{\D}{\item[({\it d})]}
\newcommand{\E}{\item[({\it e})]}
\newcommand{\G}{\item[({\it g})]}
\newcommand{\Š}{\`e}
\newcommand{\…}{\`a}
\newcommand{\•}{\`o}
\newcommand{\—}{\`u}
\newcommand{\}{\`{\i}}
\def\tag{\renewcommand{\theequation}}
\newcommand{\Footnote}{~\footnote}
\newcommand{\ie}{{\it i.e.}}
\newcommand{\dist}{\mbox{\rm dist\,}}
\newcommand{\const}{\mbox{\rm const}}
\newcommand{\trace}{\mbox{\rm trace}}
\newcommand{\Bo}{\par\hfill{$\Box$}\par\noindent}
\newcommand{\Nor}[1]{\langle{#1}\rangle_q}
\newcommand{\vs}{\vspace*{.5cm}\par\noindent}
\newcommand{\Vs}{\vspace*{.6cm}\par\noindent}
\newcommand{\Vvs}{\vspace*{.7cm}\par\noindent}
\newcommand{\VVs}{\vspace*{.8cm}\par\noindent}
\newtheorem{definition}{Definition}[section]
\newcommand{\Bd}{\begin{definition}\begin{rm}}
\newcommand{\Ed}{\end{rm}\end{definition}}
\newtheorem{remark}{Remark}[section]
\newcommand{\Br}{\begin{remark}\begin{rm}}
\newcommand{\Er}{\end{rm}\end{remark}}
\newtheorem{proposition}{Proposition}[section]
\newcommand{\Bp}{\begin{proposition}\begin{sl}}
\newcommand{\EP}[1]{\end{sl}\label{proposition:#1}\end{proposition}}
\newcommand{\propref}[1]{{\rm Proposition \ref{proposition:#1}}}
\newcommand{\Bt}{\begin{theorem}\begin{sl}}
\newcommand{\Et}{\end{sl}\end{theorem}}
\newcommand{\Bl}{\begin{lemma}\begin{sl}}
\newcommand{\El}{\end{sl}\end{lemma}}
\newtheorem{theorem}{Theorem}[section]
\newtheorem{lemma}{Lemma}[section]
\newtheorem{corollary}{Corollary}[section]
\newcommand{\eqref}[1]{{\rm (\ref{eq:#1})}}
\newcommand{\Bc}{\begin{corollary}\begin{sl}}
\newcommand{\Ec}{\end{sl}\end{corollary}}
\newcommand{\ET}[1]{\end{sl}\label{theorem:#1}\end{theorem}}
\newcommand{\EDD}[1]{\end{rm}\label{definition:#1}\end{definition}}
\newcommand{\EL}[1]{\end{sl}\label{lemma:#1}\end{lemma}}
\newcommand{\theoref}[1]{{\rm Theorem \ref{theorem:#1}}}
\newcommand{\ER}[1]{\end{rm}\label{remark:#1}\end{remark}}
\newcommand{\EC}[1]{\end{sl}\label{corollary:#1}\end{corollary}}
\newcommand{\remref}[1]{{\rm Remark \ref{remark:#1}}}
\newcommand{\cororef}[1]{{\rm Corollary \ref{corollary:#1}}}
\newcommand{\lemmref}[1]{{\rm Lemma \ref{lemma:#1}}}
\newcommand{\essup}[1]{{\rm ess}\,{{\displaystyle \sup_{\hspace*{-5mm}{#1}}}}}

\renewcommand{\i}{{\rm i}}

\pagenumbering{arabic}
\newcommand{\QED}{{\par\hfill$\square$\par}}
\renewcommand{\thefootnote}{(\arabic{footnote})}
\title{On Bifurcating Time-Periodic Flow of a Navier-Stokes Liquid past a Cylinder} 
\author{ Giovanni P. Galdi 
\thanks{Department of Mechanical Engineering and Materials Science, University of Pittsburgh, PA 15261. 
Work  partially supported by NSF DMS Grant-1311983.}}
\date{}
\maketitle
\begin{abstract}
We provide general sufficient conditions for  branching out of a time-periodic family of solutions from steady-state solutions to the two-dimensional Navier-Stokes equations in the exterior of a cylinder. To this end, we first show that the problem can be formulated as a coupled elliptic-parabolic nonlinear system in appropriate function spaces. This is obtained by separating the time-independent averaged component of the velocity field from its ``purely periodic'' one. We then  prove that time-periodic bifurcation  occurs,  provided the  linearized  time-independent operator of the parabolic problem possess a simple eigenvalue that crosses the imaginary axis when the Reynolds number passes through a (suitably defined) critical value. We also show that only supercritical or subcritical bifurcation may occur. Our approach  is different and, we believe, more direct than those used by previous authors in similar, but distinct,  context. 
 \end{abstract}

\renewcommand{\theequation}{\arabic{section}.\arabic{equation}}
\setcounter{section}{0}
\section{Introduction} 
One of the most classical phenomena in  fluid dynamics is the spontaneous oscillation of the wake in the flow of a viscous liquid past a circular cylinder. More precisely, suppose that a cylinder,  $\mathscr C$, of diameter $d$ is placed with its axis ${\sf a}$ orthogonal to the flow of a viscous liquid having an upstream constant velocity $\bfv_\infty$. Let
 $\lambda:=|\bfv_\infty|/(\nu d)$ be the relevant Reynolds number of the flow, with $\nu$ kinematic viscosity of the liquid. It is then experimentally observed that there is a critical value, $\lambda_0\sim 50$, such that if $\lambda<\lambda_0$ the motion of the liquid in a region sufficiently far from the ends of $\mathscr C$ that includes  $\mathscr C$, is  planar, steady and stable, whereas as soon as $\lambda>\lambda_0$, the motion is still planar, but its regime becomes oscillatory, as evidenced by the time-periodic motion of the wake behind $\mathscr C$; \cite[Chapter 3]{Tritton}. It is worth emphasizing that the unsteadiness of the flow arises spontaneously even though the imposed condition --uniform flow at far distances-- is time-independent. 

This is a beautiful and clear example of what, in mathematical terms, is defined as {\em time-periodic bifurcation}.
The main objective of this paper is to provide a rigorous  analysis of this  interesting phenomenon. 

In this respect, we begin to recall that, from the mathematical viewpoint, this means to investigate the following set of (dimensionless) equations  
\be\ba{cc}\medskip\left.\ba{cc}\medskip
\partial_t{\bfV}+\lambda(\bfV-\bfe_1)\cdot\nabla\bfV=\Delta\bfV-\nabla P\\
\Div\bfV=0
\ea\right\}\ \ \mbox{in $\Omega\times \real$}\\
\bfV=\bfe_1\ \ \mbox{at $\partial\Omega\times \real$}\,,\ea
\eeq{0.1}
with the further condition
\be
\lim_{|x|\to\infty}\bfV(x,t)=\0\,, \ \ \mbox{ $t\in \real$}\,.
\eeq{0.2}
Here $\bfV$ and $P$ are velocity and pressure fields of the liquid, $\Omega$ is the relevant two-dimensional unbounded region of flow (the entire portion of the plane outside the normal cross-section of $\mathscr C$), and $\bfe_1$ is a unit vector   parallel to $\bfv_\infty$. It is known that, under suitable assumptions on $\lambda_0$, the above equations  possess a unique steady-state solution branch $(\bfu(\lambda),p(\lambda))$, with $\lambda$ in  a neighborhood $U(\lambda_0)$ \cite{GaHB}. 
Writing $\bfV=\bfv(x,t;\lambda)+\bfu(x;\lambda)$, $P={\sf p}(x,t;\lambda)+p(x;\lambda)$,  equations \eqref{0.1}--\eqref{0.2} become 
\be\ba{cc}\medskip\left.\ba{rl}\medskip
\partial_\tau{\bfv}\!+\!\lambda\big[(\bfv-\bfe_1)\cdot\nabla\bfv +\bfu(\lambda)\cdot\nabla\bfv+\bfv\cdot\nabla\bfu(\lambda)\big]\!=&\!\!\!\!\Delta\bfv-\nabla {\sf p}\\
\Div\bfv=&\!\!\!\!0
\ea\!\right\}\, \ \mbox{in $\Omega\times \real$}\\
\bfv=\0\ \ \mbox{at $\partial\Omega\times\real$}\,,\ea
\eeq{1.1}
with
\be
\lim_{|x|\to\infty}\bfv(x,t)=\0\,, \ \ \mbox{ $t\in \real$}\,.
\eeq{1.2}
Our bifurcation problem consists then in finding sufficient conditions for the existence of a   non-trivial family of time-periodic solutions to \eqref{1.1}--\eqref{1.2}, $(\bfv(\lambda),{\sf p}(\lambda))$, $\lambda\in U(\lambda_0)$, of period $T=T(\lambda)$ (unknown as well), such that $(\bfv(t;\lambda),\nabla{\sf p}(t;\lambda))\to (\0,\0)$ as $\lambda\to\lambda_0$.

In order to better understand the heart of the problem and the motivation behind our approach, we 
begin to observe that, formally, \eqref{1.1}--\eqref{1.2} can be thought of as a special case of evolution equations of the type
\be
\ode{u}{t}=N(\lambda,u)\,,\ \ t\in\real
\eeq{1.3}
where $N:(\lambda,u)\in U(\lambda_0) \times X\mapsto Y$ ($X,Y$ Banach spaces) is a smooth nonlinear operator 
with $N(\lambda,0)=0$, for all $\lambda\in U(\lambda_0)$. The objective is then to find a family of time-periodic solutions $u=u(t;\lambda)$ of period $T=T(\lambda)$, $\lambda\in U(\lambda_0)$, such that $u(\lambda)\to 0$ as $\lambda\to\lambda_0$. 
Let $\mathscr L=\mathscr L(\lambda)$ be the linearization of $N(\lambda,\cdot)$ around $u=0$. Following the original ideas of E.~Hopf \cite{Hopf}, one has in mind to employ the Implicit Function Theorem, so that the bifurcation problem boils down to find conditions on $\lambda_0$,  and $T(\lambda_0)$ ensuring that the operator 
$$
\ode{}{t}-\mathscr L(\lambda_0)
$$ 
is continuously invertible in a suitable class of time-periodic functions. The latter implies that, in particular, the operator $\mathscr L(\lambda_0)$ must enjoy this property as well.

The first comprehensive  investigation of time-periodic bifurcation  may be traced back to the influential work of E.~Hopf \cite{Hopf} in the case $X=Y=\real^n$.\footnote{However, see also the previous contributions of Poincar\'e \cite{po} and Andronov \& Witt \cite{AW}.} There, Hopf shows the occurrence of bifurcation under a set of conditions that can be roughly summarized as follows: (C1) $0$ is not an eigenvalue of $L_0:=\mathscr L(\lambda_0)$, which ensures, in particular,  that  $L_0$ is continuously invertible; (C2) $L_0$ possesses two and only two purely imaginary eigenvalues $\pm{\rm i}\,\omega_0 \ (\neq 0)$ that are also simple, and (C3) As $\lambda$ passes through $\lambda_0$, the eigenvalues of
$\mathscr L(\lambda)$ ``cross'' the imaginary axis with nonzero speed.   

The approach introduced by Hopf,  lends itself to a natural extension to the infinite-dimensional case, at least when the underlying function space has a Hilbert structure, and $L_0$ is the generator of an analytic semigroup, with compact resolvent.

Along these lines of thought,  Iudovich \cite{Yu}, Joseph \& Sattinger \cite{JS}, and Iooss \cite{Io}  pioneered the investigation of the occurrence of self-oscillation in a viscous flow in a {\em bounded} domain. More precisely, they   
furnished sufficient conditions,  basically of the same type as those listed above, for the existence (and uniqueness) of bifurcating time-periodic solutions from steady-state solutions to the Navier-Stokes equations.
\footnote{For further development of the theory, its generalization  and major applications to the Navier-Stokes equations, we refer to, e.g., \cite{CR,IJ,CI,Z1,K} and the reference therein.}

It is important to emphasize that the assumption that the flow occurs in a {\em bounded domain} is fundamental. 
In fact,   it ensures, among other things, that $L_0$ has a purely discrete spectrum which, in turn, implies  that  if $0$ is not an eigenvalue then $L_0$ has a bounded inverse. 

In the case of an {\em  exterior domain} (flow past an obstacle), $L_0$ assumes the following form
\be
L_0(\bfv):={\rm P}\,\big[\Delta\bfv+\lambda_0(\bfe_1\cdot\nabla\bfv-\bfu(\lambda_0)\cdot\nabla\bfv-\bfv\cdot\nabla\bfu(\lambda_0))\big]\,,
\eeq{1.4}
where ${\rm P}$ is the Helmholtz projection.\footnote{For notation, see the next section.} However, when defined in its ``natural space'', namely, the subspace, $Z^{2,q}$, $1<q<\infty$, of solenoidal functions from the Sobolev space $W^{2,q}$ with zero trace at the boundary, the operator $L_0$ in \eqref{1.4} shows a non-empty {\em essential spectrum}, and worse, $0$ is a point of this spectrum for {\em all} $\lambda_0$  \cite{Bab,FN}. As a result, the property of continuous invertibility of $L_0$ is no longer secured.  

Nonetheless, if we define $L_0$ on an appropriate {\em homogeneous} Sobolev space, $\mathcal B$,  then we can show that $L_0$ is Fredholm of index $0$ (see \cite[Theorem 3.1]{GaR}, \cite[Theorem 1.8]{GaHB}), so that, in this framework, bounded invertibility is again guaranteed by requiring that $0$ is not an eigenvalue. These observations suggest that  for flow past an obstacle, the study of \eqref{1.1}--\eqref{1.2} and the associated  time-periodic bifurcation problem   should be performed in the Banach space, $\mathcal B$, where $L_0$ enjoys the Fredholm property.  This is, in fact, the approach employed by Babenko \cite{Bab2}, successively revisited and extended in a non-trivial way by Sazonov \cite{Saz},  also along the ideas of the seminal paper \cite{Yu}.\footnote{See also  Melcher {\em et al.} \cite{Sch} for a whole-space, vorticity formulation.}  
We wish to emphasize that the methods used by these authors work in dimension $n=3$, whereas they do not admit
any sensible generalization to the case $n=2$; see, e.g., \cite[p.39]{Bab3}.
   
It is worth remarking that the space  $\mathcal B$ above is  a subclass of the space where steady-state solutions to \eqref{0.1}--\eqref{0.2} exist. In this respect, we recall that it is a standard procedure, for flow in exterior domains, to formulate time-periodic problems in  function spaces where steady-state solutions exist; see, e.g., \cite{Ma,MaPa,Koz,Ya,GaSo,KaTs,THN} and the reference therein. However, as first pointed out in  \cite{GaD,GaA}  even though ``natural'' at first sight (``{\em steady state solution is a special case of a time-periodic one}''), this formulation is not convenient, and, in fact, as detailed in \cite{GaA},  is unable to cover the two-dimensional problem of  flow past a cylinder, which is the  focus of this paper. In view of these considerations, in \cite{GaD,GaA} we introduced a different method that, essentially, consists in reformulating the original problem as a coupled nonlinear system constituted by an elliptic equation,  formulated in the ``natural'' space of steady-state solutions,  and a parabolic equation that can be framed in a much ``better'' space; see also the analysis of Kyed in \cite{Ky}.

The approach we propose to study  time-periodic bifurcation stems from the one introduced in \cite{GaA} and, in our opinion, is simpler and more direct than those of \cite{Bab2,Saz}, with the further advantage that it allows us to cover the two-dimensional case. Basically, it consists of two main steps. We first split the sought solution, $u$, as sum of its time average over a period, $u^{(1)}$, and of its ``purely periodic'' component, $u^{(2)}$, with zero average. Accordingly,  the original problem \eqref{1.3} then splits into a coupled (nonlinear) ``elliptic-parabolic'' system in the unknowns $u^{(1)}$ and $u^{(2)}$, of the form
\be         
L_0(u^{(1)})=N_1(\lambda,u^{(1)},u^{(2)})\,,\ \ \ode{}{t}u^{(2)}-L_0(u^{(2)})=N_2(\lambda,u^{(1)},u^{(2)})\,,
\eeq{1.5} 
where $N_i$, $i=1,2$, are suitable (smooth) nonlinear operators. Now, the crucial point that makes our method different than those of \cite{Bab2,Saz}, is that, in spite of the fact that the operator $L_0$ is {\em formally} the same (see \eqref{1.4}) we frame the two equations in \eqref{1.5} in  two {\em different} function spaces, by choosing domains ${\sf D}$ and ranges ${\sf R}$ of $L_0$ appropriately. 
Specifically, in \eqref{1.5}$_1$ we take $L_0\equiv\tilde{\mathscr L}_0$ with ${\sf D}$ coinciding with the ``natural'' Banach space $\mathcal B$ of steady-state solutions and ${\sf R}\subseteq L^q$, for suitable $q>1$, whereas in \eqref{1.5}$_2$, we pick $L_0\equiv \mathscr L_0$, with ${\sf D}:= Z^{2,2}$,  and ${\sf R}\subseteq L^2$. Once the steady-state part $u^{(1)}$ of the solution has been ``isolated'' in the sense specified above, we can then show that the bifurcation problems reduces, essentially, to the study of the property of the parabolic operator $du/dt-\mathscr L_0(u)$ in \eqref{1.3}$_2$ { in the standard $L^2$ context}, exactly as in the case of a {\em bounded} region of flow \cite{Yu,JS}.  For this reason, it presents no further {\em conceptual} difficulties.          

Although our approach could be applied to a vast class of evolutionary equations,  we shall employ it here to study time-periodic bifurcation from a steady-state Navier-Stokes flow past a cylinder. More specifically, under the assumptions \eqref{H1}--\eqref{H3} formulated later on in the paper --which resemble conditions (C1)--(C3) of the original paper of Hopf--  we show, by means of the implicit function theorem, the existence of a one-parameter family of time-periodic solutions, branching out the steady-state solution ${\sf s}_0$ at $\lambda=\lambda_0$; see \theoref{3.1}(a). A characteristic feature of these solutions is that they exist either for $\lambda<\lambda_0$ or for $\lambda>\lambda_0$, so that the bifurcation is either subcritical or supercritical; see \theoref{3.1}(c). Moreover, we prove that (up to a phase shift) any other time-periodic solution branching out of $(\lambda_0,{\sf s}_0)$ must belong to the above family, under a further assumption on the branching frequency that is validated by numerical tests; see \theoref{3.1}(b) and \theoref{4.1}.

In more detail, the plan of the paper is the following. After introducing some basic notation in Section 2, in the following Section 3 our main objective is to analyze  the relevant  properties of the linearized operators $\tilde{\mathscr L}_0$ and ${\mathscr L}_0$. To this end, we recall in \propref{2.1} that $\tilde{\mathscr L}_0$, when defined in the classical  homogeneous and anisotropic Sobolev space of steady-state solutions to \eqref{1.2}--\eqref{1.2},  is Fredholm of index 0. This circumstance is supportive of our first assumption \eqref{H1}, namely, that the null space of $\tilde{\mathscr L}_0$ is trivial, so that $\tilde{\mathscr L}_0$ is boundedly invertible. Successively, we analyze the properties of the operator $\mathscr L_0$ in a subspace of the Sobolev space $W^{2,2}$, and its ``parabolic'' counterpart, $\mathscr Q(u):=du/dt-\mathscr L_0(u)$, in the space  $\mathscr W^2_{2\pi,0}$ of maximal $L^2$-regularity of $2\pi$-periodic functions with zero average over a period. In this respect, in \propref{2.3}, we prove that  $\mathscr L_0$ may have an at most countable number of purely imaginary eigenvalues that can only cluster at 0; moreover, each of these eigenvalues is isolated and of finite algebraic multiplicity. This provides the basis, on the one hand, of  our assumption \eqref{H2} that requires, in particular, that $\mathscr L_0$ has a simple, purely imaginary eigenvalue. On the other hand, by resorting to a classical result on perturbations of simple eigenvalues (see \propref{2.5}), it also supports assumption \eqref{H3} regarding the way in which the eigenvalues  of  $\mathscr L(\lambda)$ ``cross'' the imaginary axis when $\lambda$ passes $\lambda_0$.  We then study the properties of the operator $\mathscr Q$, and show that  it is Fredholm of index 0 (\lemmref{2.5}). The latter, combined with assumption \eqref{H2} allows us to give necessary and sufficient conditions for the bounded invertibility of $\mathscr Q$.  With these results in hand, in Section 4, \theoref{3.1},  under the assumptions \eqref{H1}--\eqref{H3} we prove the result of existence of a one-parameter family of time-periodic solutions to \eqref{1.1}--\eqref{1.2} mentioned earlier on, with  the property of being either subcritical or supercritical.  Finally, the uniqueness property of these solutions is discussed, in full generality, in Section 5; see \theoref{4.1}.

\section{Notation} We let $\nat$, $\mathbb Z$,  and ${\mathbb R}$, $\mathbb C$ represent, in the order,  the sets of positive and relative integers, and the fields of real and complex numbers. Thus, $\real^2$ indicates the whole plane. 
The canonical base in $\real^2$ is denoted by  $\mathfrak B:=\{\bfe_1,\bfe_2\}$. A  vector $\bfu$
will have two components in $\mathfrak B$, denoted by $u_1$ and $u_2$, respectively.
Likewise, coordinates of a point $x\in\real^2$ in the  frame $\{O,\bfe_1,\bfe_2\}$, $O\in \mathbb R^2$, will be indicated by $x_1,x_2$.

$\Omega$ stands for a fixed planar exterior domain, namely, the complement of the closure of a bounded, open,  and simply connected set, $\Omega_0$, of $\mathbb R^2$. We shall assume  $\Omega$ of class $C^2$. Moreover, we take the origin $O$ of the coordinate system in $\Omega_0$, and denote by $R_*>0$ a number such that the closure of $\Omega_0$ is strictly contained in the disk $\{\bfx\in\real^2: (x_1^2+x_2^2)^{\frac12}<R_*\}$.  

For $R\ge R_*$, we let
$$
\Omega_R=\Omega\cap \{\bfx\in\real^2: (x_1^2+x_2^2)^{\frac12}<R\}\,,\ \ \Omega^R=\Omega-\bar{\Omega_R}\,,
$$
where the bar denotes closure.

We set $\partial_t\bfu:=\partial \bfu/\partial t$, $\partial_1\bfu:=\partial \bfu/\partial x_1$, and indicate by $D^2\bfu$ the matrix of the second derivatives of $\bfu$.

For an open and connected set ${A} \subseteq {\mathbb R}^2,$ $L^q (A)$, $L^q_{loc}(A)$, $1\leq q \leq \infty,$  
$W^{m,q}({A}),$ $W_0^{m,q}(A)$, $m \geq 0,$  $(W^{0,q}\equiv W^{0,q}_0\equiv L^q$), stand for the usual Lebesgue and Sobolev classes, respectively, of real or complex functions.\Footnote{We shall use the same font style to denote scalar, vector and tensor
function spaces.}   
Norms in $L^q(A)$ and $W^{m,q}(A)$ are indicated by $\|.\|_{q,A}$ and $\|.\|_{m,q,A}$. The scalar product of functions $u,v\in L^2(A)$ will be denoted by $\langle u,v\rangle_A$. In the above notation,  the symbol $A$ will be omitted, unless confusion arises. 

As customary, for $q\in [1,\infty]$ we let $q'=q/(q-1)$ be its H\"older conjugate.  

By $D^{1,q}(\Omega)$, $1<q<\infty$, we denote the space of (equivalence classes of) functions $u$ such that
$ 
\|\nabla u\|_q<\infty\,.
$ Moreover, setting,
$$
\cald(\Omega):=\{\bfu\in C_0^\infty(\Omega):\Div\bfu=0\}
$$ 
we let $\mathcal D_0^{1,2}(\Omega)$ be the completion of $\cald(\Omega)$ in the norm $\|\nabla (\cdot)\|_2$, and set
$$
Z^{2,2}(\Omega):=W^{2,2}(\Omega)\cap \cald_0^{1,2}(\Omega)\,.
$$

Furthermore, we denote by $H_q(\Omega)$, $1<q<\infty$, ($H_2(\Omega)\equiv H(\Omega)$) the completion of $\cald(\Omega)$ in the norm $L^q(\Omega)$  and let  
 ${\rm P}_q$ be the (Helmholtz)  projection from $L^q(\Omega)$ onto $H_q(\Omega)$. ${\rm P}_q$ is independent of $q$ \cite[\S III.1]{GaB}, so that we shall simply denote it by ${\rm P}$.

For $q\in (1,3/2)$, we define
$$\ba{rl}\medskip
X^{2,q}(\Omega):=\Big\{\bfu\in L^1_{loc}(\real^2)\!\!\!&\!\!: u_2\in L^{\frac{2q}{2-q}}(\Omega)\,, \nabla u_2, \partial_1\bfu, D^2\bfu \in L^q(\Omega)\,,\\
\!\!\!&\bfu\in L^{\frac{3q}{3-2q}}(\Omega)\,,\nabla \bfu\in L^{\frac{3q}{3-q}}(\Omega)
\Big\}\,.  
\ea
$$
and 
$$
X^{2,q}_0(\Omega):=\Big\{\bfu\in X^{2,q}(\Omega): \Div\bfu=0\,,\ \bfu|_{\partial\Omega}=\0\Big\}\,.
$$
As is known, $X^{2,q}(\Omega)$ and $X^{2,q}_0(\Omega)$ become Banach spaces when endowed
with the ``natural'' norm
$$
\|\bfu\|_{X^{2,q}}:=\|u_2\|_{\frac{2q}{2-q}}+\|\nabla u_2\|_q+\|u_1\|_{\frac{3q}{3-2q}}+\|\nabla \bfu\|_{\frac{3q}{3-q}}+\|\partial_1\bfu\|_q+\|D^2\bfu\|_q\,;
$$
see \cite[\S XII.5]{GaB}.

If $M$ is a map between two spaces, we denote by ${\sf N}\,[M]$ and ${\sf R}\,[M]$ its null space and range, respectively.

In the following, $B$ is a real Banach space with associated norm $\|\cdot\|_B$.

By $B_{\mathbb C}:=B+{\rm i}\, B$ we denote the complexification of $B$.

For $q\in [1,\infty]$,  $L^q(0,2\pi;B)$ is the space of functions
$u:(0,2\pi)\rightarrow B$ such that 
$$
\left( \Int0{2\pi}\|u(t)\|_B^q \right)^{\frac 1q}<\infty, \ \ \mbox{if 
$q\in [1,\infty)\,;$}\ \ \  
\essup{t\in[0,2\pi]}\|u(t)\|_B <\infty, \ \ \mbox{if $q=\infty.$}
$$
Also, by $C(0,2\pi;B)$ we indicate the set of functions 
$u:[0,2\pi]\rightarrow B$ which are continuous from $[0,2\pi]$ with values in $B.$

Given a function $u:(x,t)\in \Omega\times[0,2\pi]\to \real^2$,
we let $\bar u=\bar u(x)$ be its average over $[0,2\pi]$, namely,
$$
{\bar u}(x):=\Frac{1}{2\pi}\int_0^{2\pi}u(x,t)dt\,.
$$
Furthermore, we shall say that
$u$ is {\em $2\pi$-periodic}, if $u(x,0)=u(x,2\pi)$, for a.a. $x\in \Omega$. Clearly, such functions  can be extended periodically to all $t\in\real$. 
We then define
$$\ba{rl}\smallskip
\mathscr W^{2}_{2\pi,0}(\Omega):=\Big\{\bfu\in
W^{1,2}(0,2\pi;H(\Omega))&\!\!\!\!\cap L^2(0,2\pi;Z^{2,2}(\Omega)):\\ &\!\!\!\!\bfu \ \mbox{is $2\pi$-periodic with $\bar\bfu=\0$}\Big\} 
\ea$$ 
with associated norm
$$
\|u\|_{\mathscr W^{2}_{2\pi,0}}:=\left(\int_0^{2\pi}\|\partial_t\bfu(t)\|_{2}^2dt\right)^{1/2}+\left(\int_0^{2\pi}\|\bfu(t)\|_{2,2}^2dt\right)^{1/2}\,.
$$
Likewise, setting $$\Omega_{2\pi}:=\Omega\times [0,2\pi]$$
we define
$$
\mathscr L_{2\pi,0}(\Omega):=\Big\{\bfu\in
L^2(\Omega_{2\pi})):\ \bfu \ \mbox{is $2\pi$-periodic with $\bar\bfu=\0$}\Big\}\,, 
$$
and its subspace
$$
\mathscr H_{2\pi,0}(\Omega):=\Big\{\bfu\in
L^2(0,{2\pi}; H(\Omega)):\ \bfu \ \mbox{is $2\pi$-periodic with $\bar\bfu=\0$}\Big\}\,. 
$$
Moreover, for $\bfu,\bfv\in \mathscr L^2_{2\pi,0}(\Omega)$ we put
$$
(\bfu|\bfv):=\int_0^{2\pi}\langle\bfu(t),\bfv(t)\rangle\,dt\,.
$$

Finally, by  $c$, $c_0$, $c_1$, etc.,  we denote positive constants, whose particular value is unessential to the context. When we wish to emphasize
the dependence of $c$ on some parameter $\xi$, we shall write  $c(\xi)$.

\setcounter{equation}{0}
\section{Analysis of the Relevant Linearized Operators}
Objective of this section is to introduce some relevant linear operators and  study their main properties in different function spaces.
 
To this end, for $\lambda_0>0$ and a given (sufficiently smooth) vector field $\bfu_0=\bfu_0(x)$, we consider the following operator, which can be viewed as a ``perturbation" to the classical Oseen operator:
\be
\bfv \mapsto \textrm{P}\,\big[\Delta\bfv+\lambda_0 \big(\partial_1{\bfv}-\bfu_0\cdot\nabla\bfv+\bfv\cdot\nabla\bfu_0\big)\big] \,.
\eeq{2.2}
The next result is shown in \cite[Theorem 1.8]{GaHB}.
\Bp Let $\bfu_0\in X^{2,q}(\Omega)$, $q\in (1,6/5]$. Then \be\tilde{\mathscr L}_0:\bfv\in X^{2,q}_0(\Omega)\mapsto 
\textrm{P}\,\big[\Delta\bfv+\lambda_0 \big(\partial_1{\bfv}-\bfu_0\cdot\nabla\bfv+\bfv\cdot\nabla\bfu_0\big)\big]\in 
H_q(\Omega)\eeq{2.1_0} is Fredholm of index 0.
\EP{2.1}

With the help of this proposition one can show the following one, whose proof can be found in \cite[Theorem 2.3]{GaHB}.
\Bp Assume that $(\bfu_0,p_0)\in X^{2,q}(\Omega)\times D^{1,q}(\Omega)$, $1<q<6/5$, is a solution to the steady-state  problem
\be 
\ba{ll}\medskip
\left.\ba{rl}\medskip
\Delta\bfu+\lambda\,\partial_1\bfu&\!\!\!=\lambda\,\bfu\cdot\nabla\bfu +\nabla p\\
\Div\bfu&\!\!\!=0\ea\right\}\ \ \mbox{in $\Omega$}\\
\bfu=\bfe_1\ \ \mbox{at $\partial\Omega$}\,,\ \ \Lim{|x|\to\infty}\bfu(x)=\0\,,
\ea
\eeq{2.2_1}
with $\lambda=\lambda_0$. Then, if
\be
{\sf N}[\tilde{\mathscr L}_0]=\{0\}\,,
\eeq{2.3}
problem \eqref{2.2_1} has a solution that is (real) analytic at $\lambda=\lambda_0$. Precisely,
there is a neighborhood $U(\lambda_0)$ of $\lambda_0$ and a family of solutions to \eqref{2.2_1},   $(\bfu(\lambda),p(\lambda))\in X^{2,q}(\Omega)\times D^{1,q}(\Omega)$,  $\lambda\in U(\lambda_0)$, such that the series
$$
\bfu(\lambda)=\bfu_0+\sum_{k=1}^\infty (\lambda-\lambda_0)^k\bfu_k\,,\ \ p(\lambda)=p_0+\sum_{k=1}^\infty (\lambda-\lambda_0)^kp_k
$$
are absolutely convergent in $X^{2,q}(\Omega)$ and $D^{1,q}(\Omega)$, respectively.
\EP{2.2}

We now consider the
operator \eqref{2.2} with domain of definition $Z^{2,2}(\Omega)\subset H(\Omega)$ and values in $H(\Omega)$, and  denote it by $\mathscr L_0$. Since $Z^{2,2}(\Omega)$ is  dense in $H(\Omega)$, $\mathscr L_0$ is densely defined.
We are interested to establish certain important properties of the spectrum $\sigma(\mathscr L_0)$. To do this, we extend $\mathscr L_0$
to a linear operator (still denoted by $\mathscr L_0$) on  $Z_{\mathbb C}^{2,2}(\Omega)$ and $H_{\mathbb C}(\Omega)$:
\be\mathscr L_0:{\sf D}_{\mathbb C}(\mathscr L_0)\subset H_{\mathbb C}(\Omega)\mapsto H_{\mathbb C}(\Omega)\,,\  \ \ {\sf D}_{\mathbb C}(\mathscr L_0):=Z^{2,2}_{\mathbb C}(\Omega)\,.\eeq{op} 
We shall then show, in particular, that the intersection of $\sigma(\mathscr L_0)$,  with $\{{\rm i}\real-\{0\}\}$ can only be constituted by a finite or countable number of eigenvalues with finite multiplicity (see \propref{2.3}).

The proof of this property requires a number of  preparatory results.
\Bl Let $\omega\in \real-\{0\}$. Then, for a given $\bff\in L^2_{\mathbb C}(\Omega)$ there is a unique corresponding $(\bfu,p)\in W_{\mathbb C}^{2,2}(\Omega)\times D^{1,2}_{\mathbb C}(\Omega)$ such that
\be 
\ba{cc}\medskip
\left.\ba{rl}\medskip
\Delta\bfu+\lambda_0\,\partial_1\bfu-\i\,\omega\,\bfu&\!\!\!=\bff +\nabla p\\
\Div\bfu&\!\!\!=0\ea\right\}\ \ \mbox{in $\Omega$}\,,\\
\bfu=\0\ \ \mbox{at $\partial\Omega$}\,.
\ea
\eeq{2.5}
Moreover, there are constants $c$ and $c_0$ depending only on $\Omega$, such that $(\bfu,p)$ satisfies the following inequality
\be
\|D^2\bfu\|_2+|\omega|^{\frac12}\|\nabla\bfu\|_2+|\omega|\|\bfu\|_2+\|\nabla p\|_2\le c\,\|\bff\|_2\,,\ \ |\omega|\ge \max\{\lambda_0^2,1\}\,.
\eeq{2.6}

\EL{2.1}
{\em Proof.} 
If we dot-multiply both sides of \eqref{2.5}$_1$ by $\bfu^*$ (${}^*:=$ complex conjugation),  integrate by parts over $\Omega$ and use \eqref{2.5}$_{2,3}$, we formally obtain
$$  
-\|\nabla\bfu\|_2^2-\i\,\omega\|\bfu\|_2^2=\langle\bff,\bfu^*\rangle-\lambda_0\langle\partial_1\bfu,\bfu^*\rangle\,. 
$$
By separating real and imaginary parts, and applying Schwartz inequality we infer
\be\ba{rl}\medskip
\|\nabla\bfu\|_2^2&\!\!\!
\le
\|\bfu\|_2\|\bff\|_2
\\
|\omega|\,\|\bfu\|_2
&\!\!\!\le\lambda_0\|\nabla\bfu\|_2+\|\bff\|_2
\ea
\eeq{2.7}
Replacing \eqref{2.7}$_1$ into \eqref{2.7}$_2$ and using Cauchy-Schwartz inequality, we easily show that
$$
|\omega|\,\|\bfu\|_2\le
\|\bff\|_2+\lambda_0\,\|\bfu\|_2^{\frac12}\|\bff\|_2^{\frac12}
\le
\left(1+\frac{\lambda_0^2}{2|\omega|}\right)\|\bff\|_2+\frac{|\omega|}{2}\,\|\bfu\|_2
$$
which implies
\be
|\omega|\,\|\bfu\|_2\le \left(2+\frac{\lambda_0^2}{|\omega|}\right)\|\bff\|_2\,.
\eeq{NaDn}
Replacing this time \eqref{NaDn} into \eqref{2.7}$_1$ we also show
\be
|\omega|^{\frac12}\|\nabla\bfu\|_2\le 
\left(2+\frac{\lambda_0^2}{|\omega|}\right)^{\frac12}\|\bff\|_2\,.
\eeq{NaDn1}
By means of the latter two estimates in conjunction with the classical Galerkin method, one can prove by standard arguments  the existence of a (weak) solution to \eqref{2.5} $(\bfu,p)\in W_{\mathbb C}^{1,2}(\Omega)\times L^2_{{\rm loc}}(\Omega_R)$ for all $R>R_*$; see, e.g. \cite[\S VII.2]{GaB}. We now write \eqref{2.5} as the following Stokes problem
\be 
\ba{cc}\medskip
\left.\ba{rl}\medskip
\Delta\bfu&\!\!\!=\bfF +\nabla p\\
\Div\bfu&\!\!\!=0\ea\right\}\ \ \mbox{in $\Omega$}\,,\\
\bfu=\0\ \ \mbox{at $\partial\Omega$}\,.
\ea
\eeq{2.8}
with $\bfF:=\lambda_0\,\partial_1\bfu-\i\,\omega\,\bfu+\bff$. In view of \eqref{NaDn}--\eqref{NaDn1} we get 
$\bfF\in L^2_{\mathbb C}$,
so that, by well-known  results  \cite[Theorem IV.5.1, Theorem V.5.3(ii)]{GaB} we deduce  that $(\bfu,p)\in  Z_{\mathbb C}^{2,2}(\Omega)\times D_{\mathbb C}^{1,2}(\Omega)$.
The existence property is thus secured. As for uniqueness  in the class $Z_{\mathbb C}^{2,2}$, it is readily established.  In fact, it is enough to proceed as in the proof of \eqref{NaDn}--\eqref{NaDn1} with $\bff\equiv\0$.
It remains to show  the validity of \eqref{2.6}. To this end, we observe that, since, obviously
$$
\|\bfF\|_2\le \big(\lambda_0\,\|\nabla\bfu\|_2+|\omega|\,\|\bfu\|_2+\|\bff\|_2\big)
\,, 
$$
from  \cite[Remark IV.4.2, Lemma  V.4.3]{GaB}
 we deduce
$$
\|\nabla p\|_2+\|D^2\bfu\|_2\le c\,\big(\lambda_0\,\|\nabla\bfu\|_2+|\omega|\,\|\bfu\|_2+\|\bff\|_2+\|\bfu\|_{2,\Omega_{R_*}}\big)\,,
$$ 
with some $c=c(\Omega)$.
Inequality \eqref{2.6} then follows, under the stated assumptions on $|\omega|$, from the latter and \eqref{NaDn}--\eqref{NaDn1}.  The   lemma is completely proved.
\QED

\Bl  The operator
$$
\mathscr K:\bfv\in Z^{2,2}(\Omega)\mapsto \bfu_0\cdot\nabla\bfv+\bfv\cdot\nabla\bfu_0\in L^2(\Omega)
$$
is compact.
\EL{2.2}
{\em Proof.}  
We begin to observe that the embeddings
\be\left.\ba{ll}\medskip
Z^{2,2}(\Omega)\subset W^{1,2}(\Omega_R)\\ 
W^{1,2}(\Omega_R)\subset L^{r}(\Omega)\,,\ \ r\in(1,\infty)
\ea\right\} \ \   \mbox{are compact, for all $R>R_*$}\,.
\eeq{2.9}
Now, let $\{\bfv_n\}\subset Z^{2,2}(\Omega)$ with $\|\bfv_n\|_{2,2}=1$, for all $n\in\nat$, and let $\bar\bfv\in Z^{2,2}(\Omega)$ be its weak limit. Without loss of generality, we may assume $\bar\bfv=\0$, which gives $\mathscr K(\bar\bfv)=\0$.
For any $R>R_*$ we show that
\be
\|\bfu_0\cdot\nabla\bfv_n\|_2\le
\|\bfu_0\|_\infty\|\nabla\bfv_n\|_{2,\Omega_R}+\|\bfu_0\|_{\infty,\Omega^R}\|\bfv_n\|_{2,2}\\
\eeq{2.10}
Likewise, by H\"older inequality,
\be
\|\bfv_n\cdot\nabla\bfu_0\|_2\le
\|\nabla\bfu_0\|_{\frac{2q}{2-q}}\|\bfv_n\|_{q^\prime,\Omega_R}+c_1\,\|\nabla\bfu_0\|_{\frac{2q}{2-q},\Omega^R}\|\bfv_n\|_{2,2}\,.
\eeq{2.11}
Since, by assumption, $\bfu_0\in X^{2,q}(\Omega)$, $1<q<6/5$,  it follows that $\bfu\in D^{1,2q/(2-q)}(\Omega)$, on the one hand, and, on the other hand, $\bfu_0\in L^\infty(\Omega)$ with $\bfu_0(x)\to \0$ uniformly, as $|x|\to \infty$; see \cite[Lemma 1]{GaA}. As a result by  \eqref{2.9}--\eqref{2.11}, and taking $R$ arbitrarily large, we may conclude
$$
\lim_{n\to\infty}\|\mathscr K(\bfv_n)\|_2=0\,.
$$
which proves
the claimed compactness property of $\mathscr K$,  and completes the proof of the proposition.\QED

\Bl Let $\bfu_0\in X^{2,q}(\Omega)$, $1<q<6/5$,  
and let $\omega\in\real-\{0\}$. Then,\footnote{By $I$ we mean the identity operator in $ H_{\mathbb C}$.} the operator
\be 
\mathscr L_0-\i\,\omega I\,,
\eeq{2.13}
with $\mathscr L_0$ defined in \eqref{op}\,,
is Fredholm of index 0.
\EL{2.3}
{\em Proof.}
We begin to notice that, as  immediately checked,
$\mathscr L_0$ is (graph) closed. In fact, this follows from \cite[Theorem 1.11 in Chapter IV]{Kato}, since $\mathscr L_0=\mathscr L_1+\mathscr K$, where $\mathscr L_1:Z^{2,2}_{\mathbb C}(\Omega)\subset H_{\mathbb C}(\Omega)\mapsto H_{\mathbb C}(\Omega)$ is obviously closed (\lemmref{2.1}) and, by \lemmref{2.2}, $\mathscr K$ is $\mathscr L_1$-compact. These two combined properties also show that \eqref{2.13} is Fredholm of index 0 (e.g. \cite[Theorem XVII.4.3]{GoGoKa}). The lemma is thus proved.\QED  

We are now in a position to show the first main result of this section.
\Bp Let $\bfu_0\in X^{2,q}(\Omega)$, $1<q<6/5$, and $\mathscr L_0$ be defined in \eqref{op}. Then $\sigma(\mathscr L_0)\cap\big\{{\rm i}\,\real-\{0\}\big\}$ consists, at
most, of a finite or countable number of eigenvalues, each of which is isolated and of finite
(algebraic)  multiplicity, that can only accumulate at 0. 
\EP{2.3}
{\em Proof.} Set $\mathscr L_\omega:=\mathscr L_0-{\rm i}\,\omega\,I$. By \lemmref{2.3} we know that $\mathscr L_\omega:H_{\mathbb C}(\Omega)\mapsto H_{\mathbb C}(\Omega)$ is an (unbounded) Fredholm operator of index 0, for all $\omega\in \real-\{0\}$. Thus, in view of well-known results (e.g. \cite[Theorem XVII.2.1]{GoGoKa}), in order to prove the stated property it is enough to show that there is $\bar\omega>0$  such that for all $|\omega|>\bar\omega$, ${\sf N}\,[\mathscr L_\omega]=\{0\}$. Now, the equation $\mathscr L_\omega(\bfv)=\0$ is equivalent to the following problem
\be
\ba{cc}\medskip
\left.\ba{rl}\medskip
\Delta\bfv+\lambda_0\,\partial_1\bfv-\i\,\omega\,\bfv&\!\!\!=\lambda_0\left(\bfu_0\cdot\nabla\bfv+\bfv\cdot\nabla\bfu_0\right) +\nabla p\\
\Div\bfv&\!\!\!=0\ea\right\}\ \ \mbox{in $\Omega$}\,,\\
\bfv=\0\ \ \mbox{at $\partial\Omega$}\,,
\ea
\eeq{2.12}
with $(\bfv,p)\in Z_{\mathbb C}^{2,2}(\Omega)\times D_{\mathbb C}^{1,2}(\Omega)$. Using \lemmref{2.1} and \eqref{2.6} in problem \eqref{2.12}, with the help of H\"older inequality we get, in particular, for all $|\omega|\ge \max\{\lambda_0^2,1\}$,
$$\ba{rl}\medskip
\|D^2\bfv\|_2+|\omega|^{\frac12}\|\nabla\bfv\|_2+|\omega|\|\bfv\|_2&\!\!\!\le c\,\lambda_0\,\|\bfu_0\cdot\nabla\bfv+\bfv\cdot\nabla\bfu_0\|_2\\ 
&\!\!\!\le c\lambda_0\left(\|\bfu_0\|_{\infty}\|\nabla\bfv\|_2+\|\nabla\bfu_0\|_{\frac{2q}{2-q}}\|\bfv\|_{q^\prime}\right)\,.
\ea$$
Using the Sobolev embedding $W^{1,2}(\Omega)\subset L^{q^\prime}(\Omega)$ in the latter, we thus infer that
$$
\|D^2\bfv\|_2+|\omega|^{\frac12}\|\nabla\bfv\|_2+|\omega|\|\bfv\|_2\le 
m_1\,\|\bfv\|_2+m_2\,\|\nabla\bfv\|_2
$$
where
$$
m_1:=c_1\,\lambda_0\,\|\nabla\bfu_0\|_{\frac{2q}{2-q}}\,,\ \ m_2=c_1\,\lambda_0\big(\|\bfu_0\|_{\infty}+\|\nabla\bfu_0\|_{\frac{2q}{2-q}}\big)
\,,
$$
and $c_1=c_1(\Omega)$,
from which the desired property follows by choosing $\bar\omega:=\max\{m_1,m_2^2,\lambda_0^2,1\}$.
\QED

Denote by $U(\lambda_0)$ a neighborhood of $\lambda_0$, and, for $q\in (1,6/5)$, let \be\lambda\in   U(\lambda_0)\mapsto \bfu(\lambda)\in X^{2,q}(\Omega)\,,\eeq{2.30}be a  continuous map  with $\bfu(\lambda_0)=\bfu_0$. Consider, alongside, the one-parameter family of operators defined by 
\be\ba{rl}\medskip
\mathscr L(\lambda):\bfv \in {\sf D}&\!\!\!\!(\mathscr L(\lambda))\subset H_{\mathbb C}(\Omega)\\
&\!\!\!\!\mapsto {\rm P}\,\big[\Delta\bfv+\lambda\,\big(\partial_1\bfv-\bfu(\lambda)\cdot\nabla\bfv-\bfv\cdot\nabla\bfu(\lambda)\big)\big]\in H_{\mathbb C}(\Omega)\,,
\ea
\eeq{2.31}
with ${\sf D}(\mathscr L(\lambda))\equiv {\sf D}(\mathscr L_0)=Z^{2,2}_{\mathbb C}(\Omega)$, and $\lambda\in U(\lambda_0)$. Obviously, 
 $\mathscr L(\lambda_0)=\mathscr L_0$.
Assume, next, that $$\mu_0:={\rm i}\,\omega_0\,,\ \ \mbox{some $\omega_0\in\real-\{0\}$\,,}$$ is in the spectrum of $\mathscr L(\lambda_0)$. Then, by \propref{2.3}, $\mu_0$ must be an eigenvalue of finite multiplicity. We are interested in the behavior of the eigenvalues, $\mu=\mu(\lambda)$, of $\mathscr L(\lambda)$ for $\lambda\in U(\lambda_0)$. To this end,
we recall that $\mu_0$ is simple (of multiplicity 1, that is) if, denoting by $\bfv_0$ the corresponding normalized eigenvector, $\bfv_0\not\in {\sf R}[\mathscr L_0-\mu_0\,I]$. Since $\mathscr L_0-\mu_0\,I$ is Fredholm of index 0, we have ${\rm dim}\, {\sf N}\,[\mathscr L_0-\mu_0\,I]={\rm codim}\,{\sf R}\,[\mathscr L_0-\mu_0\,I]=1$, and   this implies, in particular, that, letting $\mathscr L^*_0$ be the adjoint operator of $\mathscr L_0$,  
${\rm dim} \,{\sf N}[\mathscr L^*_0-\mu_0\,I]=1$ and that there is $\bfv_0^*\in {\sf N}[\mathscr L^*_0-\mu_0\,I]$ such that $\langle \bfv_0^*,\bfv_0\rangle\neq 0$; see, e.g., \cite[Section 8.4]{Z}. For convenience, we normalize $\bfv_0^*$ in such a way that
\be
\langle \bfv_0^*,\bfv_0\rangle =\pi^{-1}\,.
\eeq{2.16_0}

The following result holds (see \cite[Proposition 79.15 and Corollary 79.16]{Z1}).
\Bp Let $\mu_0$ be a simple eigenvalue of $\mathscr L_0$, and let the map \eqref{2.30} be of class $C^k$, $k\ge1$. Then, there are neighborhoods $U_1(\lambda_0)\subseteq U(\lambda_0)$ of $\lambda_0$, and $V(\mu_0)\subset \mathbb C$ of $\mu_0$, such that for each $\lambda\in U_1(\lambda_0)$ there is one and only one eigenvalue $\mu(\lambda)\in V(\mu_0)$ of $\mathscr L(\lambda)$. Moreover, the map $\lambda\mapsto\mu(\lambda)$ is of class $C^k$ and we have
\be
\mu^\prime(\lambda_0)=\langle\bfv_0^*,\partial_1\bfv_0-\bfu_0\cdot\nabla\bfv_0-\bfv_0\cdot\nabla\bfu_0-\lambda_0\,\big(\bfu^\prime(\lambda_0)\cdot\nabla\bfv_0+\bfv_0\cdot\nabla\bfu^\prime(\lambda_0)\big)\rangle\,.
\eeq{mul}
\EP{2.5} 

We now turn our focus to the study of  some important properties of the {\em time-dependent} operator
\be
\mathscr Q:=\omega_0\,\partial_{\tau}-\mathscr L_0:\ \mathscr W_{2\pi,0}^{2}(\Omega)\mapsto \mathscr H_{2\pi,0}(\Omega)\,, 
 \ \ \omega_0>0\,.
\eeq{2.17}
In particular, we are interested in determining necessary and sufficient conditions under which $\mathscr Q$ possesses a bounded inverse. 
To this end, we begin to recall the following result, which is a particular case of that proved  in \cite[Proposition 3]{GaA} 
\Bl The operator
$$
\omega_0\,\partial_{\tau}-{\rm P}\,[\Delta+\lambda_0\,\partial_1]:  \mathscr W_{2\pi,0}^{2}(\Omega)\mapsto \mathscr H_{2\pi,0}(\Omega) 
$$
is a homeomorphism.
\EL{2.4}

With the help of this lemma, we can prove the following one.
\Bl Let $\bfu_0\in X^{2,q}(\Omega)$. Then, the operator $\mathscr Q$ defined in \eqref{2.17}
is Fredholm of index 0.
\EL{2.5}
{\em Proof.} In view of \lemmref{2.4}, it is enough to show that the operator
$$
\mathscr C:\bfv\in \mathscr W_{2\pi,0}^{2}(\Omega)\mapsto \bfu_0\cdot\nabla\bfv+\bfv\cdot\nabla\bfu_0\in \mathscr L_{2\pi,0}^2(\Omega)
$$
is compact. Let $\{\bfv_k\}\subset \mathscr W_{2\pi,0}^{2}(\Omega)$ with $\|\bfv_k\|_{\mathscr W_{2\pi,0}^{2}}=1$, for all $k\in\nat$. We may then select a sequence (again denoted by $\{\bfv_k\}$) and find $\bfv_*\in \mathscr W_{2\pi,0}^{2}(\Omega)$ such that
\be
\bfv_k\to{\bfv_*} \ \ \mbox{weakly in $\mathscr W_{2\pi,0}^{2}(\Omega)$.}
\eeq{2.18}
Without loss of generality, we may take $\bfv_*\equiv\0$. From \eqref{2.18} and Lions-Aubin lemma we then have
\be
\int_0^{2\pi}\left(\|\bfv_k(\tau)\|_{2,\Omega_R}^2+\|\nabla\bfv_k(\tau)\|_{2,\Omega_R}^2\right)\to 0\ \ \mbox{as $k\to\infty$, for all $R>R_*$\,,}
\eeq{2.19}
which implies, by embedding,
\be
\int_0^{2\pi}\|\bfv_k(\tau)\|_{q^\prime,\Omega_R}^2\to 0\ \ \mbox{as $k\to\infty$, for all $R>R_*$\,,}
\eeq{2.20}
By the H\"older inequality, 
$$
\int_0^{2\pi} \|\bfu_0\cdot\nabla\bfv_k(\tau)\|_{2}^2
\le \|\bfu_0\|_\infty^2\int_0^{2\pi} \|\nabla\bfv_k(\tau)\|_{2,\Omega_R}^2+ \|\bfu_0\|_{\infty,\Omega^R}^2\int_0^{2\pi}\|\nabla\bfv_k(\tau)\|_{2}^2\,,
$$
which, by \eqref{2.19}  and the arbitrariness of $R$ implies
\be
\lim_{k\to\infty}\int_0^{2\pi} \|\bfu_0\cdot\nabla\bfv_k(\tau)\|_{2}^2=0\,.
\eeq{2.21}
Likewise, again by H\"older inequality,
$$\ba{rl}\medskip
\Int0{2\pi} \|\bfv_k(\tau)\cdot\nabla\bfu_0\|_{2}^2
\le &\!\!\!\!\|\nabla\bfu_0\|_{\frac{2q}{2-q}}^2\Int0{2\pi} \|\bfv_k(\tau)\|_{q^\prime,\Omega_R}^2\\ &+ \|\nabla\bfu_0\|_{\frac{2q}{2-q},\Omega^R}^2\Int0{2\pi}\|\bfv_k(\tau)\|_{q^\prime}^2\,.
\ea
$$
Recalling that $\mathscr W_{2\pi,0}^2(\Omega)\subset L^\infty(0,2\pi;L^{s}(\Omega))$, for all $s\in [2,\infty)$ (e.g. \cite[Lemma 2(a)]{GaA}, from the latter inequality and \eqref{2.20} we deduce
\be
\lim_{k\to\infty}\int_0^{2\pi} \|\bfv_k(\tau)\cdot\nabla\bfu_0\|_{2}^2=0\,.
\eeq{2.22}
Combining \eqref{2.21} and \eqref{2.22} we thus conclude
$$
\lim_{k\to\infty}\|\mathscr C(\bfv_k)\|_{L^2(\Omega_{2\pi})}=0\,,
$$
which completes the proof of the lemma.\QED

\Bl Let $\mathscr L_0$ be as in \propref{2.3}. Assume  $\mu_0:={\rm i}\,\omega_0\in \sigma(\mathscr L_0)$ is a simple eigenvalue, while $\mu_k:={\rm i}\,k\,\omega_0\not\in \sigma(\mathscr L_0)$, whenever $k\in\nat-\{0,1\}$.\footnote{Notice that, by \propref{2.3}, there could be only a {\em finite number} of such $\mu_k$.} Let $\bfv_0$ be the (unique) normalized eigenvector corresponding to $\mu_0$\,, and set 
\be
\bfv_1=\Re[\bfv_0\,{\rm e}^{{\rm i}\,\tau}]\,,\ \ \bfv_2=\Im[\bfv_0\,{\rm e}^{{\rm i}\,\tau}]\,. 
\eeq{2.25_0}
Then 
$$
{\rm dim}\,{\sf N}\,[\mathscr Q]=2
$$
and $\{\bfv_1,\bfv_2\}$ is a basis in $\,{\sf N}\,[\mathscr Q]$.
\EL{2.6}
{\em Proof.} It is clear that $\mathcal S:={\rm span}\,\{\bfv_1,\bfv_2\}\subseteq{\sf N}\,[\mathscr Q]$. Conversely, take $\bfw\in {\sf N}[\mathscr Q]
$. We may  expand $\bfw$ in Fourier series
\be
\bfw=\sum_{\ell=-\infty}^{\infty}\bfw_\ell\,{\rm e}^{{\rm i}\,\ell\,\tau}\,;\ \ \bfw_{\ell}(x):=\Frac{1}{2\pi}\int_{0}^{2\pi}  \bfw(x,t)\,{\rm e}^{-{\rm i}\,{\ell}\,\tau}dt
\,;\ \ \bfw_0(x)\equiv\0\,.
\eeq{Fi}
Evidently, $\bfw_\ell\in Z^{2,2}_{\mathbb C}(\Omega)\equiv{\sf D}_{\mathbb C}\,(\mathscr L_0)$\,.
From $\mathscr Q(\bfw)=\0$ we then deduce
$$
{\rm i}\,\ell\,\omega_0\,\bfw_\ell-\mathscr L_0(\bfw_\ell)=\0\,, \ \bfw_\ell\in {\sf D}_{\mathbb C}\,(\mathscr L_0)\,,\ \ \mbox{$\ell\in\mathbb Z$}\,,
$$
which, recalling that from \propref{2.3}  the $\mu_k$'s can only be eigenvalues,  by assumption and \eqref{Fi}$_3$ implies $\bfw_\ell=\0$ for all $\ell\in\mathbb Z-\{\pm 1\}$. Thus $\bfw\in \mathcal S$, and the lemma is completely proved.\QED  

We are now in a position to show the second main result of this section.
\Bp Let the assumptions of \lemmref{2.6} be satisfied, and set
\be
\bfv_1^*=\Re[\bfv_0^*\,{\rm e}^{-{\rm i}\,\tau}]\,,\ \ \bfv_2^*=\Im[\bfv_0^*\,{\rm e}^{-{\rm i}\,\tau}]\,. 
\eeq{Fra} 
where $\bfv^*_0$ is the (uniquely determined) element of ${\sf N}\,[\mathscr L_0^*-\mu_0\,I]$ satisfying \eqref{2.16_0}.\footnote{Recall that ${\rm dim}\,{\sf N}\,[\mathscr L_0^*-\mu_0\,I]=1$.} Then, for a given $\bff\in \mathscr H_{2\pi,0}(\Omega)$, necessary and sufficient condition for the problem
$$
\mathscr Q(\bfv):=\omega_0\,\partial_\tau{\bfv}-\mathscr L_0(\bfv)=\bff\,,\ \  \bfv\in \mathscr W_{2\pi,0}^2(\Omega)\,,
$$ 
to have a solution is that 
\be
(\bfv_1^*|\bff)=(\bfv_2^*|\bff)=0\,. 
\eeq{Sie}
This solution is also unique, provided 
$(\bfv_1^*|\bfv)=(\bfv_2^*|\bfv)=0$ and, in such a case, there is $c=c(\Omega)$ such that
$$
\|\bfv\|_{\mathscr W_{2\pi,0}^2}\le c\,\|\bff\|_{\mathscr H_{2\pi,0}}\,.
$$
\EP{2.5} 
{\em Proof.} Since, by \lemmref{2.5}, $\mathscr Q$ is Fredholm of index 0, and, by   \lemmref{2.6} ${\rm dim}\,{\sf N}\,[\mathscr Q]=2$, it follows (e.g. \cite[Proposition 8.14(4)]{Z}) that ${\rm dim}\,{\sf N}\,[\mathscr Q^*]=2$ where
$$
\mathscr Q^*:=-\omega_0\,\partial_{\tau}-\mathscr L_0^*
$$ 
is the adjoint of $\mathscr Q$. In view of the stated properties of $\bfv_0^*$\,, we infer that ${\rm span}\,\{\bfv_1^*,\bfv_2^*\}= {\sf N}\,[\mathscr Q^*]$, and the proposition  follows from another  classical result on Fredholm operators (e.g., \cite[Proposition 8.14(2)]{Z}).\QED
\section{Bifurcating Time-Periodic Solutions.}
We begin to put the original problem \eqref{1.1}--\eqref{1.2} in a different and equivalent form that will allow us to employ the results established in the previous section. 

To this end, let $\lambda_0$ $(>0)$ be a value of the Reynolds number for which the steady-state problem  \eqref{2.3} has a solution $(\bfu_0,p_0)\in X^{2,q}(\Omega)\times D^{1,q}(\Omega)$, $1<q<6/5$. We suppose that $\lambda_0$ is such that $(\bfu_0,p_0)$ is a point of  an analytic solution branch $(\bfu(\lambda),p(\lambda))$ to \eqref{2.3}, for all $\lambda$ in a neighborhood $U(\lambda_0)$ of $\lambda_0$.
By \propref{2.2}, such a $\lambda_0$ exists if we assume that
\tag{H1}\be
\mbox{${\sf N}\,[\tilde{\mathscr L}_0]=\{0\}$, with $\tilde{\mathscr L}_0$ defined in  \eqref{2.1_0}}\,, 
\eeq{H1}
or, {\em equivalently},
\tag{H1$^*$}
\be
\ba{ll}\medskip
\left.\ba{rl}\medskip
\Delta\bfu+\lambda_0\big(\,\partial_1\bfv&\!\!\!- \bfu_0\cdot\nabla\bfu -\bfu\cdot\nabla\bfu_0\big)=\nabla \phi\\
\Div\bfu&\!\!\!=0\ea\right\}\,\ \mbox{in $\Omega$}\\
\bfu=\0 \ \mbox{at $\partial\Omega$}\,,\ \ (\bfu,\phi)\in X^{2,q}(\Omega)\times D^{1,q}(\Omega)\,,\ 1<q<\frac65,
\ea  \Longrightarrow \ \bfu\equiv \0\,.
\eeq{H1*}
\setcounter{equation}{0}\renewcommand{\theequation}{\arabic{section}.\arabic{equation}}\par Our first objective is to prove the existence of a family of time-periodic solutions of  period $T:=2\pi/\omega$ (to be determined) to \eqref{1.1}, bifurcating from the point $(\lambda_0; (\bfu_0,p_0))$.
To this end, following  Lindtstedt \cite{L} and  Poincar\'e \cite{po}, we introduce  the scaled time $\tau:=\omega\,t$, so that \eqref{1.1} becomes
\be\ba{cc}\medskip\left.\ba{rl}\medskip
\omega\,\partial_\tau{\bfv}\!+\!\lambda\big[(\bfv-\bfe_1)\cdot\nabla\bfv +\bfu(\lambda)\cdot\nabla\bfv+\bfv\cdot\nabla\bfu(\lambda)\big]\!=&\!\!\!\!\Delta\bfv-\nabla {\sf p}\\
\Div\bfv=&\!\!\!\!0
\ea\!\right\}\, \ \mbox{in $\Omega_{2\pi}$}\\
\bfv=\0\ \ \mbox{at $\partial\Omega_{2\pi}$}\,,\ea
\eeq{3.1}
 
We next split $\bfv$ and ${\sf p}$ as the sum of their time average, $(\bar{\bfv},\bar {\sf p})$, over the time interval $[0,2/\pi]$, and their ``purely periodic" component $(\bfw:=\bfv-\bar{\bfv},\varphi:=\bar{\sf p}-{\sf p})$. In this way, problem \eqref{3.1} can be  equivalently rewritten as the following coupled nonlinear elliptic-parabolic problem
\be 
\ba{c}\medskip
\left.\ba{rl}\medskip
\Delta\bar{\bfv}+\lambda_0\big(\,\partial_1\bar\bfv-\bfu_0\cdot\nabla\bar\bfv-\bfu_0\cdot\nabla\bar\bfv\big)&\!\!\!= \nabla \bar{\sf p}+\bfN_1(\lambda,\bar{\bfv},\bfw)\\
\Div\bar\bfv&\!\!\!=0\ea\right\}\ \ \mbox{in $\Omega$}\\
\bar\bfv=\0\ \ \mbox{at $\partial\Omega$}\,,
\ea 
\eeq{3.2}
and
\be\ba{cc}\medskip\left.\ba{rl}\smallskip
\omega_0\,\partial_\tau{\bfw}-\Delta\bfw-\lambda_0\,\big(\partial_1\bfw-&\!\!\!\!\bfu_0\cdot\nabla\bfw-\bfw\cdot\nabla\bfu_0\big)\\ \medskip
&\!\!\!\!=\nabla {\varphi} +\bfN_2(\lambda,\omega,\bar\bfv,\bfw)\\
\Div\bfw&\!\!\!\!=0
\ea\right\}\ \ \mbox{in $\Omega_{2\pi}$}\\
\quad\quad\bfw=\0\ \ \mbox{at $\partial\Omega_{2\pi}$}\,,\ea
\eeq{3.3}
where $\omega_0>0$, and
\be\ba{rl}\medskip
\bfN_1:=&\!\!\!\!(\lambda_0-\lambda)\,\left[\partial_1\bar\bfv-\bfu(\lambda)\cdot\nabla\bar\bfv-\bar\bfv\cdot\nabla\bfu(\lambda)\right]\\ \medskip
&\!\!\!\!+\lambda_0\,\big[(\bfu(\lambda)-\bfu_0)\cdot\nabla\bar\bfv+\bar\bfv\cdot\nabla(\bfu(\lambda)-\bfu_0)\big]\\
&\!\!\!\!+\lambda\Big[\bar\bfv\cdot\nabla\bar\bfv+
\bar{\bfw\cdot\nabla\bfw}\,\Big]\,,
\ea
\eeq{3.4}
and
\be\ba{rl}\medskip
\bfN_2:=&\!\!\!\! (\omega_0-\omega)\,\partial_\tau\bfw +(\lambda-\lambda_0)\,\left[\partial_1\bfw-\bfu(\lambda)\cdot\nabla\bfw-\bfw\cdot\nabla\bfu(\lambda)\right]\\ \smallskip
&\!\!\!\!-\lambda_0\,\Big[(\bfu(\lambda)-\bfu_0)\cdot\nabla\bfw+\bfw\cdot\nabla(\bfu(\lambda)-\bfu_0)\Big]
\\
&\!\!\!\!+\lambda\Big[\bfw\cdot\nabla\bar\bfv+\bar\bfv\cdot\nabla\bfw+\bfw\cdot\nabla\bfw-
\bar{\bfw\cdot\nabla\bfw}\,\Big]
\ea
\eeq{3.5}

Some functional properties of the quantities $\bfN_i$, $i=1,2$, are proved next.
\Bl Let $1<q<6/5$. The following bilinear maps are continuous 
$$\ba{ll}\medskip
\calm_1 :(\bfv_1,\bfv_2)\in [X^{2,q}(\Omega)]^2\mapsto \bfv_1\cdot\nabla\bfv_2\in L^q(\Omega)\,,\\ \medskip 
\calm_2: (\bfw_1,\bfw_2)\in[\mathscr W_{2\pi,0}^2(\Omega)]^2\mapsto 
\Int0{2\pi}
{\bfw_1\cdot\nabla\bfw_2}\in  L^r(\Omega)\,,\ r=q,2\,,\\ \medskip
\calm_3:(\bfv,\bfw)\in X^{2,q}(\Omega)\times \mathscr W_{2\pi,0}^2(\Omega)\mapsto \bfv\cdot\nabla\bfw\in \mathscr L^2_{2\pi,0}(\Omega)\,,\\ \medskip
\calm_4:(\bfv,\bfw)\in X^{2,q}(\Omega)\times \mathscr W_{2\pi,0}^2(\Omega)\mapsto \bfw\cdot\nabla\bfv\in \mathscr L^2_{2\pi,0}(\Omega)\,,\\
\calm_5:(\bfw_1,\bfw_2)\in [\mathscr W_{2\pi,0}^2(\Omega)]^2\mapsto \bfw_1\cdot\nabla\bfw_2\in \mathscr L^2_{2\pi,0}(\Omega)\,.
\ea
$$
\EL{3.1}
{\em Proof.} The continuity of $\calm_1$ is shown in \cite[Lemma XII.5.4]{GaB}. 
In order to show the remaining properties, we recall the following continuous embeddings (see \cite[Lemmas 1 and 2]{GaA})
\be\ba{cc}\medskip
X^{2,q}(\Omega)\subset L^\infty(\Omega)\,;\ \ X^{2,q}(\Omega)\subset D^{1,\frac{2q}{2-q}}(\Omega)\,;\\
\mathscr W_{2\pi,0}^2(\Omega)\!\subset\! L^\infty(0,2\pi; L^s(\Omega))\,, \mbox{all $s\in [2,\infty)$}\,;\, \ \mathscr W_{2\pi,0}^2(\Omega)\!\subset\! L^4(0,2\pi; W^{1,4}(\Omega))\,. 
\ea
\eeq{3.6}
Therefore, by \eqref{3.6} and  H\"older inequality we deduce
\be\ba{ll}\medskip
\|\calm_2(\bfw_1,\bfw_2)\|_q\le \Int{0}{2\pi}\|\bfw_1\|_{\frac{2q}{2-q}}\|\nabla\bfw_2\|_2\le c_1\,\|\bfw_1\|_{\mathscr W_{2\pi,0}^2}\|\bfw_2\|_{\mathscr W_{2\pi,0}^2}\\ \medskip
\|\calm_2(\bfw_1,\bfw_2)\|_2\le \Int{0}{2\pi}\|\bfw_1\|_{4}\|\nabla\bfw_2\|_4\le c_2\,\|\bfw_1\|_{\mathscr W_{2\pi,0}^2}\|\bfw_2\|_{\mathscr W_{2\pi,0}^2}\\  \medskip
\|\calm_3(\bfw,\bfw)\|_{\mathscr L^2_{2\pi,0}}\le \|\bfv\|_\infty\Big(\Int{0}{2\pi}\|\nabla\bfw_2\|_2^2\Big)^{\frac12}\le c_3\,\|\bfv\|_{X^{2,q}}\|\bfw_2\|_{\mathscr W_{2\pi,0}^2}\\
\medskip
\|\calm_4(\bfw,\bfv)\|_2\le \|\nabla\bfv\|_{\frac{2q}{2-q}}\Big(\Int{0}{2\pi}\|\bfw\|_{q^\prime}^2\Big)^{\frac12}\le c_4\,\|\bfv\|_{X^{2,q}}\|\bfw\|_{\mathscr W_{2\pi,0}^2}\\
\|\calm_5(\bfw_1,\bfw_2)\|_2\le \Big(\Int{0}{2\pi}\|\bfw_1\|_{4}^4\Big)^{\frac14}\Big(\Int{0}{2\pi}\|\nabla\bfw_2\|_{4}^4\Big)^{\frac14}\!\!\le c_5\,\|\bfw_1\|_{\mathscr W_{2\pi,0}^2}\|\bfw_2\|_{\mathscr W_{2\pi,0}^2}
\ea
\eeq{ineq}
\QED

With the help of this lemma, and recalling the definition of $\tilde{\mathscr L}_0$ and $\mathscr Q$ given in \eqref{2.3} and \eqref{2.19}, respectively, we may infer that problems \eqref{3.2}--\eqref{3.5} can be {\em equivalently} rewritten in the following operator form
\be\ba{ll}\medskip 
\tilde{\mathscr L}_0(\bar\bfv)=\mathcal N_1(\lambda,\bar\bfv,\bfw)\ \ \mbox{in $H_q(\Omega)$}\,,\\
\mathscr Q(\bfw)=\mathcal N_2(\lambda,\omega,\bar\bfv,\bfw)\ \ \mbox{in $\mathscr H_{2\pi,0}(\Omega)$}\,,
\ea
\eeq{3.7}
where $\mathcal N_i={\rm P}\,\!\bfN_i$, $i=1,2$. 

The  desired bifurcation result will be obtained by showing that, under appropriate assumptions on $(\lambda_0,\omega_0)$,  there exists  a non-trivial family of solutions $(\bar\bfv, \bfw)\in X^{2,q}_0(\Omega)\times \mathscr W^2_{2\pi,0}(\Omega)$ to \eqref{3.7} for  $(\lambda,\omega)$ in a neighborhood of $(\lambda_0,\omega_0)$. 
\Br The asymptotic side condition \eqref{1.2} is embodied in the function spaces where $\bar\bfv$ and $\bfw$ are sought. In fact, since $\bar\bfv\in X^{2,q}_0(\Omega)$, from \cite[Lemma 1]{GaA} we have
$$
\lim_{|x|\to\infty}|\bar\bfv(x)|=0 \ \ \mbox{uniformly},
$$
whereas $\bfw\in\mathscr W_{2\pi,0}^2(\Omega)$ and \cite[Theorem II.9.1]{GaB} imply, for almost all $t\in [0,2\pi]$,
$$
\lim_{|x|\to\infty}|\bfw(x,t)|=0 \ \ \mbox{uniformly in $x$}.
$$
\ER{4.1}

We next recall that, by \propref{2.3}, if ${\rm i}\,\omega_0\in \sigma(\mathscr L_0)$  then it must be an eigenvalue of finite multiplicity, and, moreover, there is at most,  a finite number of eigenvalues of the form ${\rm i}\,k\,\omega_0$ with $k\in\nat$. With this in mind, we make the more stringent hypothesis that $\omega_0$ be such that
\tag{H2}
\be\ba{ll}\medskip
\mu_0:={\rm i}\,\omega_0 \,\ \mbox{is  an eigenvalue of multiplicity 1 of $\mathscr L_0$}\,,\\ 
k\,\mu_0\,, k\in\nat-\{0,1\}\ \mbox{is not an eigenvalue of $\mathscr L_0$}\,,
\ea
\eeq{H2}
and look for solutions to \eqref{3.7} satisfying the further requirement
\renewcommand{\theequation}{\arabic{section}.\arabic{equation}}\setcounter{equation}{7}
\be
(\bfw|\bfv_1^*)=\varepsilon\,,\ \ (\bfw|\bfv_2^*)=0\,, 
\eeq{3.8}
where $\bfv_i^*$, $i=1,2$, is defined in \propref{2.5}, and $\varepsilon\in (-1,1)$.
\par 
We are now in a position to prove our main result on the existence and  uniqueness  of bifurcating time-periodic solutions, along with their relevant properties. To this end, we observe that, under the assumptions \eqref{H1} and \eqref{H2} the eigenvalue $\mu(\lambda)$ of the operator $\mathscr L(\lambda)$ defined in \eqref{2.17} is a  $C^\infty$-function of $\lambda$ in a suitable neighborhood of $\lambda_0$, and \eqref{mul} holds.
\Bt Suppose that \eqref{H1} and \eqref{H2} hold and that, in addition,
\tag{H3}\be
\Re[\mu'(\lambda_0)]\neq 0\,.
\eeq{H3}
\renewcommand{\theequation}{\arabic{section}.\arabic{equation}}\setcounter{equation}{11}Then,  the following properties are valid. 
\begin{itemize}
\item[{\rm (a)}] {\rm Existence.}  
There are (real) analytic families 
\be
\big(\bar\bfv(\varepsilon),\bfw(\varepsilon),\omega(\varepsilon),\lambda(\varepsilon)\big)\in X^{2,q}_0(\Omega)\times \mathscr W_{2\pi,0}^2(\Omega)\times \real_+^2
\eeq{fam}
satisfying \eqref{3.7}--\eqref{3.8} for all $\varepsilon$ in a neighborhood $\mathcal I(0)$ of $0$, and such that
\be
(\bar\bfv(\varepsilon),\bfw(\varepsilon)-\varepsilon\,\bfv_1,\omega(\varepsilon),\lambda(\varepsilon))\to
(\0,\0, \omega_0,\lambda_0) \ \ 
\mbox{as $\varepsilon\to 0$} 
\eeq{dis} 
with $\bfv_1$ given in \eqref{2.25_0}\,. Moreover, the corresponding  velocity field $\bfV$ of  the original problem \eqref{1.1} has the following form near $\varepsilon=0$
\be 
\bfV(x,\tau;\lambda(\varepsilon))=\bfu_0(\bfx)+\varepsilon\,\big[(\cos \tau)\bfa_1
+(\sin  \tau)\bfa_2\big]+\varepsilon^2\,[\bfV_1+\bfV_2]\,,
\eeq{3.12} 
where $\bfa_i\in Z^{2,2}(\Omega)$, $i=1,2$, and $(\bfV_1,\bfV_2)\in X^{2,q}_0(\Omega)\times\mathscr W^2_{2\pi,0}(\Omega)$ satisfy
$$
\|\bfV_1\|_{X^{2q}}+\|\bfV_2\|_{\mathscr W^2_{2\pi,0}}\le M\,,
$$
with $M$  independent of $\varepsilon\to 0$.
\item[{\rm (b)}] {\rm Uniqueness.} There is a neighborhood  $$\mathscr U(\0,\0,\omega_0,\lambda_0)\subset X^{2,q}_0(\Omega)\times \mathscr W_{2\pi,0}^2(\Omega)\times \real_+^2$$such that every (nontrivial) solution to \eqref{3.7} lying in $\mathscr U$ must belong, up to a phase shift, to the family \eqref{fam}.
\item[{\rm (c)}] {\rm Parity.} The functions $\omega(\varepsilon)$ and $\lambda(\varepsilon)$ are even:
$$
\omega(\varepsilon)=\omega(-\varepsilon)\,,\ \ \lambda(\varepsilon)=\lambda(-\varepsilon)\,,\ \ \mbox{for all $\varepsilon\in\cali(0)$\,.} 
$$
Consequently, the bifurcation due to these solutions is either subcritical or supercritical, a two-sided bifurcation being excluded.\footnote{Unless $\lambda\equiv \lambda_0$.}
\end{itemize}

\ET{3.1}
{\em Proof.} We rescale our original unknowns in \eqref{3.7}--\eqref{3.8} as follows:
\be
\bfw=\varepsilon\, {\sf w}\,,\ \ \bar\bfv=\varepsilon\,{\sf v}\,, 
\eeq{stpo}
so that  \eqref{3.7}--\eqref{3.8}
can be  equivalently written as
\newcommand{\sfv}{{\sf v}} 
\newcommand{\sfw}{{\sf w}}
\be
\ba{lcc}\medskip 
\tilde{\mathscr L}_0(\sfv)-\mathscr N_1(\varepsilon,\lambda,\sfv,\sfw)=0\ \ \mbox{in $H_q(\Omega)$}\,,\\ \medskip
\mathscr Q(\sfw)-\mathscr N_2(\varepsilon,\lambda,\omega,\sfv,\sfw)=0\ \ \mbox{in $\mathscr L^2_{2\pi,0}(\Omega)$}\,,
\\ ({\sf w}|\bfv_1^*)=1\,,\ \ ({\sf w}|\bfv_2^*)=0\,,\ea
\eeq{3.9} 
where
\be\ba{rl}\medskip
\mathscr N_1:=&\!\!\!\!{\rm P}\,\Big\{(\lambda_0-\lambda)\,\left[\partial_1\sfv-\bfu(\lambda)\cdot\nabla\sfv-\sfv\cdot\nabla\bfu(\lambda)\right]\\ \medskip
&\!\!\!\!+\lambda_0\,\big[(\bfu(\lambda)-\bfu_0)\cdot\nabla\sfv+\sfv\cdot\nabla(\bfu(\lambda)-\bfu_0)\big]\\
&\!\!\!\!+\lambda\,\varepsilon\,\big[\sfv\cdot\nabla\sfv+
\bar{\sfw\cdot\nabla\sfw}\,\big]\Big\}\,,
\ea
\eeq{3.10}
and
\be\ba{rl}\medskip
\mathscr N_2:=&\!\!\!\! {\rm P}\,\Big\{(\omega_0-\omega)\,\partial_t\sfw +(\lambda-\lambda_0)\,\left[\partial_1\sfw-\bfu(\lambda)\cdot\nabla\sfw-\sfw\cdot\nabla\bfu(\lambda)\right]\\ \smallskip
&\!\!\!\!-\lambda_0\,\Big[(\bfu(\lambda)-\bfu_0)\cdot\nabla\sfw+\sfw\cdot\nabla(\bfu(\lambda)-\bfu_0)\Big]
\\
&\!\!\!\!+\lambda\,\varepsilon\,\big[\sfw\cdot\nabla\sfv+\sfv\cdot\nabla\sfw+\sfw\cdot\nabla\sfw-
\bar{\sfw\cdot\nabla\sfw}\,\big]\Big\}
\ea
\eeq{3.11}
Define the map:
$$\ba{cc}\medskip
F: (\varepsilon,\lambda, \omega,\sfv,\sfw)\in \cali(0)\times U(\lambda_0)\times V(\omega_0)\times X^{2,q}_0(\Omega)\times \mathscr W_{2\pi,0}^2(\Omega)\\ \smallskip
\mapsto
\Big(\tilde{\mathscr L}_0(\sfv)-\mathscr N_1(\varepsilon,\lambda,\sfv,\sfw),\ 
\mathscr Q(\sfw)-\mathscr N_2(\varepsilon,\lambda,\omega,\sfv,\sfw),\
({\sf w}|\bfv_1^*)-1,\ ({\sf w}|\bfv_2^*)\Big)\\
\in H_q(\Omega)\times \mathscr H_{2\pi,0}(\Omega)\times \real^2\,.
\ea
$$
The nonlinear terms in \eqref{3.10}--\eqref{3.11} are of polynomial form, and, by \eqref{H1} and \propref{2.2}, $\bfu(\lambda)$ is analytic. Thus, also with the help of \lemmref{3.1} we may conclude that  $F$ is analytic. Furthermore, from \eqref{3.9}--\eqref{3.11} and \eqref{H1} it follows that for $\varepsilon =0$,  the equation $F=0$ has the solution $(\lambda=\lambda_0,\omega=\omega_0,\sfv=\0,\sfw=\bfv_1)$. Therefore, by the (real) analytic version of the implicit function theorem (e.g. \cite[Proposition 8.11]{Z}), to show the existence part in the theorem -including the validity of \eqref{dis}- it suffices to show that the Fr\'echet derivative of $F$ with respect to ${\sf U}:=(\lambda,\omega,\sfv,\sfw)$ evaluated at $(\varepsilon=0,\lambda=\lambda_0,\omega=\omega_0,\sfv=\0,\sfw=\bfv_1$) is a bijection. The latter will hold if  we  prove that for any
$({\sf f}_1,{\sf f}_2,{\sf f}_3,{\sf f}_4)\in H_q(\Omega)\times\mathscr H_{2,\pi,0}(\Omega)\times\real\times\real$, the following set of equations has one and only one solution $(\lambda,\omega,\sfv,\sfw)\in \real\times\real\times X^{2,q}_0(\Omega)\times \mathscr W_{2\pi,0}^2(\Omega)$:
\be\ba{rl}\medskip
\tilde{\mathscr L}_0(\sfv)=&\!\!\!\!{\sf f}_1\ \ \mbox{in $H_{q}(\Omega)$}\\ \medskip
\mathscr Q(\sfw)= &\!\!\!\!\mathscr F(\lambda,\omega,\bfv_1)+{\sf f}_2\ \ \mbox{in $\mathscr H_{2\pi,0}(\Omega)$}\,,\\
({\sf w}|\bfv_1^*)=&\!\!\!\!{\sf f}_3\,,\ \ ({\sf w}|\bfv_2^*)={\sf f}_4 \ \ \mbox{in $\real$}\,,
\ea
\eeq{3.13}
where
\be\ba{rl}\medskip
\mathscr F(\lambda,\omega,\bfv_1):=-&\!\!\!\!\omega\,\partial_\tau\bfv_1+\lambda\,{\rm P}\,\Big\{\partial_1\bfv_1-\bfu_0\cdot\nabla\bfv_1-\bfv_1\cdot\nabla\bfu_0
\\
&\!\!\!\! \quad\quad
-\lambda_0\big(\bfu'(\lambda_0)\cdot\nabla\bfv_1+\bfv_1\cdot\nabla\bfu'(\lambda_0)\big)\Big\}\,.\ea
\eeq{3.13_1}
Since $\tilde{\mathscr L}_0$ is Fredholm of index 0 (\propref{2.1}), in view of \eqref{H1}, for any given ${\sf f}_1\in H_q(\Omega)$, equation \eqref{3.13}$_1$ has one and only one solution $\sfv\in X^{2,q}(\Omega)$. Therefore, it remains to prove the existence and uniqueness property only for the system of equations \eqref{3.13}$_{2-4}$
To this aim, we observe that, by \propref{2.5}, for a given ${\sf f}_2\in \mathscr H_{2\pi,0}(\Omega)$,  equation \eqref{3.13}$_2$ possesses a solution $\sfw\in\mathscr W_{2\pi,0}^2(\Omega)$ if and only if its right-hand side satisfies \eqref{Sie}. By a direct calculation, from \eqref{2.25_0},  \eqref{Fra}, and \eqref{2.16_0} we show
\be
(\partial_\tau\bfv_1|\bfv_1^*)=0\,,\ \ (\partial_\tau\bfv_1|\bfv_2^*)=\pi\,. 
\eeq{3.14}
Furthermore, again by a straightforward calculation that uses also \eqref{mul} and the fact that ${\rm P}$ is self-adjoint in $L^2$ and ${\rm P}\,\bfv_1^*=\bfv_1^*$, we infer
\be\ba{rr}\medskip
\left\langle{\rm P}\,\big\{\partial_1\bfv_1-\bfu_0\cdot\nabla\bfv_1-\bfv_1\cdot\nabla\bfu_0
+\lambda_0\big(\bfu'(\lambda_0)\cdot\nabla\bfv_1+\bfv_1\cdot\nabla\bfu'(\lambda_0)\big)\big\},\bfv_1^*\right\rangle\\
=\Re[\mu'(\lambda_0)]\,.
\ea
\eeq{3.15}
Employing \eqref{3.14} and \eqref{3.15}, we thus recognize that the compatibility condition \eqref{Sie} for the solvability of equation \eqref{3.13}$_2$ reduces to
solving the following algebraic system for  $\lambda$ and $\omega$:
\be\ba{rl}\medskip
\lambda\,\Re[\mu'(\lambda_0)]&\!\!\!\!=-({\sf f}_2|\bfv_1^*)\\
-\omega\,\pi+\lambda\, (\mathscr F|\bfv_2^*)&\!\!\!\!=-({\sf f}_2|\bfv_2^*)\,,
\ea
\eeq{3.16}
with $\mathscr F$ given in \eqref{3.13_1}. By virtue of \eqref{H3}, for any given ${\sf f}_2$ in the specified class, we can always find (uniquely determined) $\lambda$ and $\omega$ satisfying the above system, and this, by \propref{2.5}, ensures the existence of a solution ${\sf w}_1$ to \eqref{3.13}$_2$ corresponding to the  selected values of $\lambda$ and $\omega$. We now set
$$ 
\sfw:={\sf w}_1+\alpha\,\bfv_1+\beta\,\bfv_2\,,\ \ \alpha\,,\, \beta\in\real\,.
$$
Clearly, by \lemmref{2.6}, ${\sfw}$ is also a solution to \eqref{3.13}$_2$. We then 
choose $\alpha$ and $\beta$ in such a way that $\sfw$ satisfies both conditions \eqref{3.13}$_{3,4}$  for any given ${\sf f}_i\in\real$, $i=1,2$. This choice is made possible by the fact that, as is immediately checked,
\be
(\bfv_i|\bfv_j^*)=\delta_{ij}\,,\ \ i,j=1,2\,.
\eeq{3.17}
The existence part is therefore accomplished. We now turn to uniqueness and set ${\sf f}_i=\0$ in \eqref{3.13}$_{2-4}$. From \eqref{3.16} and \eqref{H3} it then follows $\lambda=\omega=0$ which in turn implies, by \eqref{3.13}$_2$ and  \lemmref{2.6}, $\sfw=\gamma_1\,\bfv_1+\gamma_2\,\bfv_2$, for some $\gamma_i\in\real$, $i=1,2$. Replacing this information back in \eqref{3.13}$_{3,4}$ with ${\sf f}_3={\sf f}_4=\0$, and using \eqref{3.17} we conclude $\gamma_1=\gamma_2=0$, which completes the uniqueness proof.
We have thus shown that the above specified Fr\'echet derivative of $F$ is a bijection, which ensures  existence to \eqref{3.9}--\eqref{3.11},  and therefore of a family of solutions, parametrized in $\varepsilon$, to \eqref{3.7}--\eqref{3.8} in the sense  specified in (a).  To complete  the proof of the statement in (a), it   remains to show \eqref{3.12}. To this end, we begin to notice that from \eqref{2.25_0} we have
\be
\bfv_1=(\cos \tau)\bfa_1+(\sin \tau)\bfa_2\,,\ \ \bfa_1,\bfa_2\in Z^{2,2}(\Omega)\,.
\eeq{v1}
Next, let us  give for granted, momentarily, the result in (c). By the analyticity property of $\lambda(\varepsilon)$ we then infer  that either $\lambda(\varepsilon)\equiv\lambda_0$ or else there is an integer $k\ge 1$ such that 
\be
\lambda(\varepsilon)=\lambda_0+\varepsilon^{2k}\lambda_k+O(\varepsilon^{2k+2})
\ \ \lambda_k
\in\real-\{0\}
\,.
\eeq{quad}
As a result, by \propref{2.2} and \eqref{H1} we deduce, in particular,
\be
\bfu(\lambda)-\bfu_0=\varepsilon^{2}\,\bfU\,,\ \ \|\bfU\|_{X^{2,q}}\le M\,,
\eeq{quad1}
with $M$ independent of $\varepsilon\to 0$. Likewise, from the analyticity properties of ${\sf w}$ and $\sfv$ and \eqref{dis} we have 
\be
\sfw-\bfv_1=\varepsilon\,{\sf W}\,,\ \ \sfv=\varepsilon\,{\sf V}\,, \ \ 
\|{\sf V}\|_{X^{2,q}}+\|{\sf W}\|_{\mathscr W_{2\pi,0}^2}^2\le M\,.
\eeq{wv}
Thus, since
$$
\bfV=\bfu_0+\varepsilon\,\bfv_1+\varepsilon\,\big[(\sfw(\varepsilon)-\bfv_1)+\sfv(\varepsilon)\big] +\bfu(\lambda)-\bfu_0\,,  
$$  
\eqref{3.12} is a consequence of this identity and \eqref{v1}--\eqref{wv}. We shall next prove the uniqueness property in (b) by adapting to the present case the abstract argument of \cite[Theorem 8.B]{Z}. 
Let $\bfz=\bar\bfz+\bfq$, $\bfq:=\bfz-\bar\bfz$  be a  $2\pi$-periodic function
where $\bar\bfz\in X^{2,q}_0(\Omega)$ and $\bfq\in \mathscr W_{2\pi,0}^2(\Omega)$ satisfy the first and the second equation in \eqref{3.7}, respectively, with $\omega\equiv\tilde{\omega}$ and $\lambda\equiv\tilde\lambda$.  
By the uniqueness property associated with  the implicit function theorem, the proof of the claimed uniqueness
amounts to show that we can find a sufficiently small $\rho>0$ such that if
\be
\|\bar\bfz\|_{X^{2,q}}+\|\bfq\|_{\mathscr W^2_{2\pi,0}}+|\tilde\omega-\omega_0|+|\tilde\lambda
-\lambda_0|<\rho\,,
\eeq{3.20}
then there exists a neighborhood of $0$, $\cali(0)\subset\real$, such that
\be\ba{cc}\medskip
\bfq=\eta\,\bfv_1+\eta\,{\sf y}\,,\, \bar\bfz=\eta\,{\sf z}\,, \ \mbox{for all $\eta\in\cali(0)$},\, \\
|\tilde\omega-\omega_0|+|\tilde\lambda-\lambda_0|+\|{\sf z}\|_{X^{2,q}}+\|{\sf y}\|_{\mathscr W^2_{2\pi,0}}\to 0\ \ \mbox{as $\eta\to 0$}\,.\ea
\eeq{3.21}
To this end, we notice that, by \eqref{3.17}, we may write
\be
\bfq=\textsf{\bf v} +\bfy
\eeq{3.22}
where ${\bf v}=(\bfy|\bfv^*_1)\,\bfv_1+(\bfy|\bfv^*_2)\,\bfv_2$ and 
\be
(\bfy|\bfv^*_i)=0\,,\ \ i=1,2\,.
\eeq{3.22_1}
We next make the simple but important observation that if we modify $\bfq$ by a  constant phase shift in time, $\delta$, namely, $\bfq(\tau)\to\bfq(\tau+\delta)$ it follows that the shifted function is still  a $2\pi$-periodic solution to \eqref{3.7}$_2$ and, moreover, by an appropriate choice of $\delta$, 
\be   
{\bf v}=\eta\,\bfv_1\,,
\eeq{3.23}
with $\eta=\eta(\delta)\in\real$. (The proof of \eqref{3.23} is straightforward, once we take into account \eqref{2.25_0}.)
Notice that from \eqref{3.20}, \eqref{3.22}--\eqref{3.23} it follows that
\be\ba{rl}
|\eta| &\!\!\!\!\le c_0\,\rho\,,\ \ c_0>0\,,\\
\|\bfy\|_{\mathscr W_{2\pi,0}^2}&\!\!\!\!\le \rho_1\,,\ \ \rho_1\to 0 \,\ \mbox{as $\rho\to 0$}\,.\ea
\eeq{3.23_1} 
From \eqref{3.7},  \eqref{3.22}, and \eqref{3.23}, we thus infer 
\be 
\tilde{\mathscr L}_0(\bar\bfz)= \mathcal N_1(\tilde\lambda,\bar{\bfz},\eta\,\bfv_1+\bfy) 
\eeq{3.24}
and
\be
\mathscr Q(\bfy)=\mathscr F\big(\eta(\tilde\omega-\omega_0),\eta(\tilde\lambda-\lambda_0),\bfv_1\big)+\mathcal N(\eta,\tilde\lambda,\tilde\omega,\bar\bfz,\bfy)
\eeq{3.25}
In \eqref{3.24} the quantities  $\caln_1$ and $\mathscr F$ are defined in  \eqref{3.7}$_1$ and \eqref{3.13_1}, respectively, whereas
$$\ba{rl}\medskip
\caln:=&\!\!\!{\rm P}\Big\{-\eta\,(\tilde\lambda-\lambda_0)\big[(\bfu(\tilde\lambda)-\bfu_0)\cdot\nabla\bfv_1+\bfv_1\cdot\nabla(\bfu(\tilde\lambda)-\bfu_0)\big]\\ \medskip
&\!\!\!-\eta\,\lambda_0\big[(\bfu(\tilde\lambda)-\bfu_0-(\tilde\lambda-\lambda_0)\bfu'(\lambda_0))\cdot\nabla\bfv_1\\ \medskip&\!\!\!+\bfv_1\cdot\nabla(\bfu(\tilde\lambda)-\bfu_0-(\tilde\lambda-\lambda_0)\bfu'(\lambda_0)\big]+\eta\,\tilde{\lambda}\,\big[\bfv_1\cdot\nabla\bfy+\bfy\cdot\nabla\bfv_1\\ \medskip
&\!\!\!+\bfv_1\cdot\nabla\bar{\bfz}+\bar{\bfz}\cdot\nabla\bfv_1 +\eta\,\bfv_1\cdot\nabla\bfv_1-\big(\eta\,\bar{\bfv_1\cdot\nabla\bfv_1}+\bar{\bfv_1\cdot\nabla\bfy}+\bar{\bfy\cdot\nabla\bfv_1}\big)\big]\Big\}\\
&\!\!\!+\caln_2(\tilde{\lambda},\tilde{\omega},\bar{\bfz},\bfy)\,,
\ea
$$
with $\caln_2$ given in \eqref{3.7}.
We now observe the following facts.\par 
(1) By \eqref{H1} and \propref{2.2},
$$\ba{rl}\medskip
\|\bfu(\tilde\lambda)-\bfu(\lambda_0)\|_{X^{2,q}}\le &\!\!\!\!M\,|\tilde\lambda-\lambda_0|
\\
\|\bfu(\tilde\lambda)-\bfu_0-(\tilde\lambda-\lambda_0)\bfu'(\lambda_0)\|_{X^{2,q}}\le &\!\!\!\!M\,|\tilde\lambda-\lambda_0|^2\,,
\ea
$$
where $M$ is independent of $|\tilde\lambda-\lambda_0|\to 0$. 
\par
(2) Since $\tilde{\mathscr L}_0$ is Fredholm of index 0, again by \eqref{H1}, it is boundedly invertible.\par
(3) By \lemmref{3.1} and (1) we easily show that
\be\ba{rl}\medskip
\|\caln_1(\tilde\lambda,\bar\bfz,&\!\!\!\!\!\eta\,\bfv_1+\bfy)\|_{H_q}\\
&\!\!\!\!\le c_1 \left(|\tilde\lambda-\lambda_0|\|\bar\bfz\|_{X^{2,q}}+\|\bar\bfz\|_{X^{2,q}}^2+\eta^2+|\eta|\,\|\bfy\|_{\mathscr W_{2\pi,0}^2}+\|\bfy\|_{\mathscr W_{2\pi,0}^2}^2\right)\,.\ea
\eeq{3.37}
and
\be\ba{rl}\medskip
\|\caln-\caln_2(\tilde\lambda,&\!\!\!\!\!\tilde\omega,\bar{\bfz},\bfy)\|_{\mathscr H_{2\pi,0}}\\
&\!\!\!\!\le c_2 \left(|\eta|\,|\tilde\lambda-\lambda_0|^2+|\eta|\,\|\bar\bfz\|_{X^{2,q}}+\eta^2+|\eta|\,\|\bfy\|_{\mathscr W_{2\pi,0}^2}\right)\,.\ea
\eeq{3.38}
Likewise,
\be\ba{rl}
\|\caln_2(\tilde\lambda,&\!\!\!\!\!\tilde\omega,\bar{\bfz},\bfy)\|_{\mathscr H_{2\pi,0}}\\ &\!\!\!\!\le c_3\,
\left(\big(|\tilde\omega-\omega_0|+|\tilde\lambda-\lambda_0|\big)\,\|\bfy\|_{\mathscr W_{2\pi,0}}+\|\bar\bfz\|_{X^{2,q}}\|\bfy\|_{\mathscr W_{2\pi,0}^2}+\|\bfy\|_{\mathscr W_{2\pi,0}^2}^2\right)\,.
\ea
\eeq{3.39}\par
(4) By \propref{2.5} and \eqref{3.15} we infer 
\be
\ba{rl}\medskip
\eta(\tilde\lambda-\lambda_0)\,\Re[\mu'(\lambda_0)]&\!\!\!\!=-({\caln }|\bfv_1^*)\\
-\eta(\tilde\omega-\omega_0)\,\pi+\eta(\tilde\lambda-\lambda_0)\, (\mathscr F|\bfv_2^*)&\!\!\!\!=-(\caln|\bfv_2^*)\,,
\ea
\eeq{basta}
where the quantity 
$$
\mathscr F=\mathscr F\big(\eta(\tilde\omega-\omega_0),\eta(\tilde\lambda-\lambda_0),\bfv_1)\,,
$$
defined in \eqref{3.13_1}, satisfies
$$
\|\mathscr F\big(\eta(\tilde\omega-\omega_0),\eta(\tilde\lambda-\lambda_0),\bfv_1\big)\|_{\mathscr H_{2\pi,0}}\le c_4\,|\eta|(|\tilde\lambda-\lambda_0|+|\tilde\omega-\omega_0|)\,.
$$
\par
(5) \propref{2.5},  \eqref{3.22_1}, and \eqref{3.25} imply that 
$$
\|\bfy\|_{\mathscr W^2_{2\pi,0}}\le c_5\left(\|\mathscr F\big(\eta(\tilde\omega-\omega_0),\eta(\tilde\lambda-\lambda_0),\bfv_1\big)+\mathcal N(\eta,\tilde\lambda,\tilde\omega,\bar\bfz,\bfy)\|_{\mathscr H_{2\pi,0}}\right)\,.
$$

With all the properties in (1)--(5) being established, we may now draw the following consequences. In the first place, by choosing $\rho$ sufficiently small and employing \eqref{3.23_1},   from \eqref{3.24} and \eqref{3.25} we deduce
\be
\|\bar\bfz\|_{X^{2,q}}+\|\bfy\|_{\mathscr W_{2\pi,0}^2}\le c_6\,\left(\eta^2+|\eta|\,(|\tilde\lambda-\lambda_0|+|\tilde\omega-\omega_0|)\right)\,.
\eeq{3.41}
Moreover, using also \eqref{H3} and \eqref{basta}, we show
\be\ba{rl}\medskip
|\eta|\,(|\tilde\lambda-\lambda_0|+|\tilde\omega-&\!\!\!\!\omega_0|)\le c_7\,\Big[|\eta|\,(|\bar\bfz\|_{X^{2,q}}+\|\bfy\|_{\mathscr W^2_{2\pi,0}})+\eta^2\\
&\!\!\!\!+(|\tilde\omega-\omega_0|+|\tilde\lambda-\lambda_0|)\|\bfy\|_{\mathscr W^2_{2\pi,0}}+\|\bar\bfz\|_{X^{2,q}}^2+\|\bfy\|_{\mathscr W^2_{2\pi,0}}^2\Big]\,.
\ea
\eeq{3.42}
Thus, combining \eqref{3.41} and \eqref{3.42}, and taking $\rho$ sufficiently small we obtain, again with the help of \eqref{3.23_1},
\be
\|\bar\bfz\|_{X^{2,q}}+\|\bfy\|_{\mathscr W_{2\pi,0}^2}\le c_8\,|\eta|^2\,,
\eeq{3.43} 
which, once used back into \eqref{3.42}, gives also
\be
|\tilde\lambda-\lambda_0|+|\tilde\omega-\omega_0|\le c_9\,|\eta|\,.
\eeq{3.44}
Recalling \eqref{3.22} and \eqref{3.23}, by virtue of \eqref{3.43} and \eqref{3.44} we may establish the validity of \eqref{3.21}, thus concluding the proof of the uniqueness property (b). It remains to show the statement in (c). To this end, we observe that if $\bfv(t):=\bar\bfv+\bfw$ is a solution to \eqref{3.7} in the function class specified in part (a),  so is $\bfv^\prime:=\bfv(t+\pi)$. Let $\bfw^\prime:=\bfv^\prime-\bar\bfv^\prime$. By the uniqueness property of part (b), we  must have (with the obvious meaning of the symbols)
$\omega^\prime(\varepsilon)=\omega(\varepsilon)$ and $\lambda^\prime(\varepsilon)=\lambda(\varepsilon)$, for all $\varepsilon $ in a neighborhood of $0$. However, if $(\bfw|\bfv_1^*)=\varepsilon\,\bfv_1$, then $(\bfw^\prime|\bfv_1^*)=-\varepsilon\,\bfv_1$.  from which the stated parity condition follows. Finally, if $\lambda\not\equiv0$, the expansion \eqref{quad} must hold, and this implies $\lambda(\varepsilon)<\lambda_0$ or $\lambda(\varepsilon)>\lambda_0$, according to whether $\lambda_k$  (the first nonzero coefficient in the Taylor expansion for $\lambda$) is negative or positive. The theorem is completely proved.
\QED
\Br
As is well known, condition \eqref{H3} means that
when the Reynolds number $\lambda$ passes through the ``critical'' value $\lambda_0$ the eigenvalues of $\mathscr L(\lambda)$  must cross the imaginary axis at $\pm {\rm i}\,\omega_0$ with  non-zero speed.
\ER{4.2}

\setcounter{equation}{0}
\section{Further Properties of Bifurcating Solutions} The results of \theoref{3.1} ensure the existence and uniqueness of bifurcating time-periodic solutions in a neighborhood, $\mathscr I$, of $(\lambda_0;(\bfu_0,p_0))$,  with period $2\pi/\omega$ and $\omega$ ``sufficiently close'' to the imaginary part, $\omega_0$ of a (simple) purely imaginary eigenvalue of the relevant linearized operator $\mathscr L_0$. These solutions are of particular physical interest in that they may branch out only sub- or super-critically. However, the theorem cannot exclude the existence of  other bifurcating time-periodic solutions in the same neighborhood $\mathscr I$, but with frequency  not ``close'' to $\omega_0$, and having, in principle, different branching properties. Objective of this section is to prove that, under suitable further assumptions, the solutions determined in \theoref{3.1} are, in fact, the {\em only} possible bifurcating time-periodic solutions in $\mathscr I$. Roughly speaking, these assumptions amount to say that as $\lambda$ passes $\lambda_0$, there are two {\em and only two} (complex conjugate) eigenvalues of the operator $\mathscr L_0$ crossing the imaginary axis, and that, in addition, for any nontrivial sequence $\{\bfv_n,\lambda_n,\omega_n\}$ of solutions to \eqref{3.1} with
$$|\bar\bfv_n\|_{X^{2,q}}+\|\bfw_n\|_{\mathscr W^2_{2\pi,0}}+|\lambda_n-\lambda_0|\to 0
$$
there exists $\delta>0$ (which may depend on the sequence) such that
\tag{H4}
\be
 \omega_n\ge \delta\,,\ \ \mbox{for all $n\in\nat$}.   
\eeq{H4}
From the physical viewpoint, this means that  time-periodic bifurcating solutions branch out with a finite (nonzero) frequency.
Numerical tests confirm that the above assumptions are indeed satisfied \cite[Section 6]{Ran},\cite{FGQ}.

\setcounter{equation}{0}\renewcommand{\theequation}{\arabic{section}.\arabic{equation}}
In order to prove the uniqueness result previously described, we need to show several preparatory lemmas, along the ideas developed in \cite{Yu}. For $\omega>0$, let
$$
\mathscr J_\omega:=\int_0^{2\pi}\left(\omega^2\,\|\partial_\tau\bfw(\tau)\|_2^2+\|{\rm P}\Delta\bfw(\tau)\|_2^2\right)\,d\tau\,.
$$
\Bl Let $\bfw\in \mathscr W_{2\pi,0}^2(\Omega)$. Then there is $c=c(\Omega)>0$ such that
\be\ba{rl}\medskip
\Int0{2\pi}\|\bfw(\tau)\|_r^2d\tau&\!\!\!\!\le c\,\Frac{1}{\omega^{\frac{2+r}{r}}}\mathscr J_\omega\,,\ \ \mbox{all $r\ge 2$}\,,\\ \medskip
\Max{\tau\in [0,2\pi]}\|\nabla\bfw(\tau)\|_2^2 
&\!\!\!\!\le \Frac{1}{\omega}\mathscr J_\omega\,,
\\
\Int0{2\pi}\|\bfw(\tau)\|_4^4d\tau&\!\!\!\!\le c\,\Frac{1}{\omega^3}\mathscr J_\omega^2\,.
\ea
\eeq{4.1} 
\EL{4.1}
{\em Proof.} We begin to observe that since $\bar\bfw(x)=\0$ for a.a. $x\in\Omega$, by the Wirtinger inequality we have
$$
\int_0^{2\pi}|\bfw(x,\tau)|^2d\tau\le \int_0^{2\pi}|\partial_\tau\bfw(x,\tau)|^2d\tau\,,
$$
which, in turn, after integration over $\Omega$ and using Fubini's theorem implies
\be
\int_0^{2\pi}\|\bfw(\tau)\|_2^2d\tau\le \frac{1}{\omega^2}\mathscr J_\omega\,. 
\eeq{4.2}
Next, we notice that from the obvious identity 
$$
\|\nabla\bfw\|_2^2=-\langle{\rm P}\Delta\bfw,\bfw\rangle
$$
and the Schwartz inequality it follows that
$$
\int_0^{2\pi}\|\nabla\bfw(\tau)\|_2^2d\tau\le\left(\int_0^{2\pi}\|{\rm P}\Delta\bfw(\tau)\|_2^2d\tau\right)^{\frac12}\left(\int_0^{2\pi}\|\bfw(\tau)\|_2^2d\tau\right)^{\frac12}\,,
$$
which in conjunction with \eqref{4.2}, by the Cauchy-Schwartz inequality delivers
\be
\int_0^{2\pi}\|\nabla\bfw(\tau)\|_2^2d\tau\le \frac{1}{2\omega}\mathscr J_\omega\,. 
\eeq{4.3}
Now, by classical interpolation, we have $\mathscr W_{2\pi,0}^2(\Omega)\subset C([0,2\pi];W^{1,2}(\Omega))$ and, for all $0\le s\le\tau\le 2\pi$,
\be
\|\nabla\bfw(\tau)\|_2^2-\|\nabla\bfw(s)\|_2^2=\int_s^\tau\langle\partial_\xi\bfw(\xi),{\rm P}\Delta\bfw(\xi)\rangle d\xi\,;
\eeq{4.4}
see \cite[Chapter 3.1]{LM}. Applying first the Schwartz inequality on the right-hand side of this equation and then integrating over $s\in [0,2\pi]$ we find for all $\tau\in [0,2\pi]$
$$ 
\|\nabla\bfw(\tau)\|_2^2\le 
\int_0^{2\pi}\|\nabla\bfw(\xi)\|_2^2d\xi+
\left(\int_0^{2\pi}\|\partial_\xi\bfw(\xi)\|_2^2d\xi\right)^{\frac12}\left(\int_0^{2\pi}\|{\rm P}\Delta\bfw(\xi)\|_2^2d\xi\right)^{\frac12}
\,.$$
Inequality \eqref{4.1}$_2$ is then a consequence of the latter and \eqref{4.3}. Furthermore, we recall the well-known embedding inequality (see, e.g., \cite[Lemma II.3.1]{GaB})
\be
\|\bfw\|_r\le c_1\|\bfw\|_2^{1-\lambda}\|\nabla\bfw\|_2^\lambda\,,\ \ \lambda:=\frac{r-2}{r}\,,\ \ r\ge 2\,.
\eeq{4.5}
Squaring both sides of the latter, integrating over $[0,2\pi]$ and using H\"older inequality, we get
$$
\int_0^{2\pi}\|\bfw(\tau)\|_r^2d\tau\le c_2
\left(\int_0^{2\pi}\|\bfw(\tau)\|_2^2d\tau\right)^{1-\lambda}\left(\int_0^{2\pi}\|\nabla\bfw(\tau)\|_2^2d\tau\right)^{\lambda}\,,
$$
which, by virtue of \eqref{4.2} and \eqref{4.3} implies \eqref{4.1}$_1$.
Finally, choosing $r=4$ in \eqref{4.5} raising both sides to the power 4 and integrating over $[0,2\pi]$ we show
$$
\int_0^{2\pi}\|\bfw(\tau)\|_4^4d\tau\le c_3
\,\max_{\tau\in[0,2\pi]}\|\bfw(\tau)\|_2^2\int_0^{2\pi}\|\bfw(\tau)\|_2^2d\tau
$$
which with the help of \eqref{4.1}$_{1,2}$ furnishes \eqref{4.1}$_3$.
\QED

\Bl Let $\omega>0$, and let  $\bfv$ be a $2\pi$-periodic solution to \eqref{3.1} with $\bar\bfv\in X^{2q}_0(\Omega)$ and $\bfw:=\bfv-\bar\bfv\in \mathscr W_{2\pi,0}^2(\Omega)$. Suppose  that
$$
\|\bar\bfv\|_{X^{2q}}+\|\bfw\|_{\mathscr W_{2\pi,0}^2}\le\rho\,,
$$
for some $\rho>0$. 
Then
$$
\omega\le \sqrt{2A+A^2}\,,
$$
where $A:=\lambda^2(C_1+C_2(\rho+\rho^2))$, with $C_1=C_1(\Omega,\|\bfu\|_{X^{2,q}})>0$ and $C_2=C_2(\Omega)>0$.
\EL{4.2}
{\em Proof.} From \eqref{3.1} we deduce that $\bfw$ satisfies the following equation
\be\ba{rl}\medskip
\omega\,\partial_\tau\bfw-{\rm P}\Delta\bfw=&\!\!\!\!\lambda\,{\rm P}\big[\partial_1\bfw-\bfu(\lambda)\cdot\nabla\bfw-\bfw\cdot\nabla\bfu(\lambda)\big]
\\ \medskip
&+\lambda\,{\rm P}\big[\bfw\cdot\nabla\bar\bfv+\bar\bfv\cdot\nabla\bfw+\bfw\cdot\nabla\bfw-\bar{\bfw\cdot\nabla\bfw}\big]\\ :=&\!\!\!\!\lambda\,(\bfg_1+\bfg_2)\,.
\ea
\eeq{4.6}
Squaring both sides of \eqref{4.6}, integrating over $[0,2\pi]$, and observing that, by the $2\pi$-periodicity and \eqref{4.4}
$$
\int_0^{2\pi}\langle\partial_\tau\bfw,\Delta\bfw\rangle d\,\tau=0\,,
$$
we infer
\be
\mathscr J_\omega=\lambda^2\int_0^{2\pi}\|\bfg_1(\tau)+\bfg_2(\tau)\|_2^2d\tau\,.
\eeq{4.7}
Arguing as in \eqref{ineq}, for any $q\in(1,6/5)$ we  show that
$$
\int_0^{2\pi}\|\bfg_1(\tau)\|_2^2d\tau\le c_1\,\int_0^{2\pi}\left(\|\nabla\bfw(\tau)\|_2^2+\|\bfw(\tau)\|_{q'}^2\right)d\tau\,,
$$
where, $c_1=c_1(\Omega,\|\bfu\|_{X^{2,q}})$. The latter inequality, in conjunction with \eqref{4.1} furnishes
\be
\int_0^{2\pi}\|\bfg_1(\tau)\|_2^2d\tau\le c_2\left(\frac{1}{\omega^{\frac{2+q'}{q'}}}+\frac{1}{\omega}\right)\mathscr J_\omega\,,
\eeq{4.8}
where, here and in the rest of the proof, $c_i$, $i=2,\ldots$,  denotes a positive constant depending at most on $\Omega$.  
Likewise, we show
\be
\int_0^{2\pi}\|\bfw\cdot\nabla\bar\bfv(\tau)+\bar\bfv\cdot\nabla\bfw(\tau)\|_2^2d\tau\le c_3\rho^2 \left(\frac{1}{\omega^{\frac{2+q'}{q'}}}+\frac{1}{\omega}\right)\mathscr J_\omega\,.
\eeq{4.9}
Moreover
$$
\int_0^{2\pi}\|\bfw\cdot\nabla\bfw(\tau)\|_2^2d\tau\le \left(\int_0^{2\pi}\|\bfw(\tau)\|_4^4d\tau\right)^{\frac12}\left(\int_0^{2\pi}\|\nabla\bfw(\tau)\|_4^4d\tau\right)^{\frac12}\,,
$$
so that by \eqref{3.6}, \eqref{4.1}$_3$ and by assumption we conclude
\be
\int_0^{2\pi}\|\bfw\cdot\nabla\bfw(\tau)\|_2^2d\tau\le c_4\,\rho\, \frac{1}{\omega^{\frac{3}{2}}}\mathscr J_\omega\,.
\eeq{4.10}
Finally, we observe that
$$
\int_0^{2\pi}\|\bar{\bfw\cdot\nabla\bfw}\|_2^2d\tau\le 2\pi \int_0^{2\pi}\|\bfw\cdot\nabla\bfw(\tau)\|_2^2d\tau\,,
$$
so that by \eqref{4.10} it follows
\be
\int_0^{2\pi}\|\bar{\bfw\cdot\nabla\bfw}\|_2^2d\tau \le c_5\,\rho\, \frac{1}{\omega^{\frac{3}{2}}}\mathscr J_\omega\,.
\eeq{4.11}
Collecting \eqref{4.8}--\eqref{4.11} and observing that, by Young's inequality,
$$
\frac{1}{\omega^{\frac{2+q'}{q'}}}=\frac{1}{\omega^{\frac{4}{q'}}}\,\frac{1}{\omega^{\frac{q'-2}{q'}}}\le c(q)\,\left(\frac{1}{\omega^{2}}+\frac{1}{\omega}\right)\,,\ \ \frac{1}{\omega^{\frac32}}=\frac{1}{\omega}\,\frac{1}{\sqrt{\omega}}\le \half\,\left(\frac{1}{\omega^{2}}+\frac{1}{\omega}\right)\,,
$$
we find
$$
\int_0^{2\pi}\|\bfg_1(\tau)+\bfg_2(\tau)\|_2^2
\le (C_1+C_7(\rho+\rho^2))\left(\frac{1}{\omega^2}+\frac{1}{\omega}\right)\mathscr J_\omega\,.
$$
Combining this latter inequality with \eqref{4.7}, we get
\be
\frac{1}{\omega^2}+\frac{1}{\omega}\ge \frac1A\,,
\eeq{4.12}
where $A$ is the quantity defined in the statement of the lemma. The proof is then accomplished by employing in \eqref{4.12}  the elementary inequality
$$
\frac1\omega\le \frac1{2A}+\frac {A^2}{2\omega^2}\,. 
$$ 
\QED
\Bl Let $\omega_*>0$ and let
\be
\mathscr Q_*:\bfw\in\mathscr W_{2\pi,0}^2(\Omega)\mapsto \omega_*\partial_\tau\bfw-\mathscr L_0(\bfw)\in\mathscr H_{2\pi,0}(\Omega)\,.
\eeq{4.13}
Then, $\mathscr Q_*$ is boundedly invertible if and only if $\mu_{*k}:={\rm i}\,k\,\omega_*\not\in\sigma(\mathscr L_0)$ for all $k\in\nat-\{0\}$. 
\EL{4.3}
{\em Proof.} From \propref{2.3} we know that the $\mu_{*k}$'s can only be eigenvalues of $\mathscr L_0$.~\footnote{In fact, the same proposition guarantees that there is at most  a {\em finite} number of such $\mu_{*k}$'s.} Since, by assumption, ${\sf N}[\mathscr L_0\pm\mu_{*k}\,I]=\{0\}$, for all $k$, with the help of \lemmref{2.3} we deduce  that the operator
$$  
(\mathscr L_0\pm\mu_{*k}\,I)^{-1} 
$$
is a homeomorphism  of $H_{\mathbb C}(\Omega)$ onto $Z_{\mathbb C}^{2,2}(\Omega)$. Therefore, by using classical Fourier series techniques, we show that $\mathscr Q_*$ is also boundedly invertible with 
$$  
\mathscr Q_*^{-1}:\bff\in \mathscr H_{2\pi,0}(\Omega)\mapsto \bfw\in\mathscr W_{2\pi,0}^2(\Omega)\,,  
$$
where
$$
\bfw(t):=\Sum{\ell=-\infty, \ell\neq0}{\infty}{\rm e}^{{\rm i}\,\ell\,\tau}\big(\mathscr L_0-{\rm i}\,{\ell}\,\omega_*\,I\big)^{-1}\big[\frac1{2\pi}\int_0^{2\pi}\bff(s)
{\rm e}^{-{\rm i}\,\ell\,s}\,ds\big]\,.
$$
\QED

We are now in a position to prove the main result of this section.
\Bt Suppose that \eqref{H1}, \eqref{H3}, and \eqref{H4} hold, and that the intersection of the spectrum of the operator $\mathscr L_0$ with the imaginary axis consists of two and only two (complex conjugate) simple eigenvalues $\pm{\rm i}\,\omega_0$. Then, there exists $\rho>0$ such that every non-trivial $2\pi$-periodic solution $\bfv$ to \eqref{3.1} for some $\omega>0$, for which
\be
\|\bar\bfv\|_{X^{2,q}}+\|\bfw\|_{\mathscr W_{2\pi,0}^2}+|\lambda-\lambda_0|<\rho\,,
\eeq{4.13}   
must belong (up to a phase shift) to the one-parameter family of solutions constructed in \theoref{3.1}. 
\ET{4.1}
{\em Proof.} In view of the uniqueness result of \theoref{3.1}(b), if the claim is not true,  there should exist a (non-trivial) sequence of solutions to \eqref{3.1}, 
$\{\bfv_n,\lambda_n,\omega_n\}$,  and a number $a>0$ such that  \be\|\bar\bfv_n\|_{X^{2,q}}+\|\bfw_n\|_{\mathscr W_{2\pi,0}^2}+|\lambda_n-\lambda_0|\to 0\,,\eeq{4.15}
and
$$
|\omega_n-\omega_0|\ge a.$$ 
By virtue of \lemmref{4.2}, the sequence $\{\omega_n\}$ is bounded, so that, by \eqref{H4}, there exists $\omega_*>0$, such that 
\be
|\omega_n-\omega_*|\to 0\,,\ \  \omega_*\neq\omega_0.
\eeq{4.16} 
From \eqref{3.1} we thus obtain that $\bar\bfv$ and $\bfw$ solve the following coupled equations
$$\ba{ll}\medskip
\tilde{\mathscr L}(\bar\bfv)=\caln_1(\lambda_0+\sigma,\bar\bfv,\bfw)\ \ \mbox{in $H_q(\Omega)$}\,,
\\
\mathscr Q_*(\bfw)=\caln_{*2}(\lambda_0+\sigma,\omega_*+\xi,\bar\bfv,\bfw)\ \ \mbox{in $\mathscr H_{2\pi,0}^2(\Omega)$}\,,
\ea
$$
where $\sigma:=\lambda-\lambda_0$, $\xi:=\omega-\omega_*$,   $\caln_1$ and $\mathscr Q_*$ are defined in  \eqref{3.7}$_1$ and \eqref{4.13}, respectively, whereas $\caln_{*2}$ is given in  \eqref{3.7}$_2$ with $\omega_0\equiv\omega_*$. Consider the map 
$$\ba{cc}\medskip
\mathcal M: (\sigma,\xi,\bar\bfv,\bfw)\in \cali(0)\times U(0)\times X^{2,q}_0(\Omega)\times \mathscr W_{2\pi,0}^2(\Omega)\\ \smallskip
\mapsto
\Big(\tilde{\mathscr L}_0(\bar\bfv)-\mathcal N_1(\lambda_0+\sigma,\bar\bfv,\bfw),\ 
\mathscr Q_*(\sfw)-\mathcal N_{*2}(\lambda_0+\sigma,\omega_*+\xi,\bar\bfv,\bfw)\Big)\\
\in H_q(\Omega)\times \mathscr H_{2\pi,0}(\Omega)\,.
\ea
$$
Clearly, the equation $\mathcal M(0,0,\bar\bfv,\bfw)=\0$ has the solution ${\sf U}_0:=(\bar\bfv=\0,\bfw=\0)$. The Fr\'echet derivative, $D\mathcal M$, of $\mathcal M$ with respect to $(\bar\bfv,\bfw)$ evaluated at $(\sigma=0,\xi=0,{\sf U}_0)$ is given by
$$
D\mathcal M: ({\sf v},{\sf w})\in X^{2,q}_0(\Omega)\times \mathscr W_{2\pi,0}^2(\Omega)\mapsto \big(\tilde{\mathscr L}_0({\sf v}),\mathscr Q_*({\sf w})\big) \in H_{q}(\Omega)\times \mathscr H_{2\pi,0}(\Omega)\,,
$$
which, by \eqref{H1} and \lemmref{4.3} is a bijection. Therefore, by the implicit function theorem, there are no nontrivial solutions satisfying \eqref{4.15} and \eqref{4.16}, thus showing a contradiction. 
\QED

\ed
By \eqref{2.25_0} we have 
$$
\bfv_1=(\cos t)\,\bfa_1+(\sin t)\,\bfa_2\,,\ \ \bfv_2=(\sin t)\,\bfa_1-(\cos t)\,\bfa_2\,,\ \ \bfa_1,\bfa_2\in Z^{2,2}(\Omega)\,.
$$
Denoting by $\omega_n$ the frequency of $\bfv_n$, we set
$$
\bfv_{1n}=(\cos \omega_nt)\,\bfa_1+(\sin \omega_nt)\,\bfa_2\,,\ \ \bfv_{2n}=(\sin \omega_nt)\,\bfa_1-(\cos \omega_nt)\,\bfa_2\,.
$$ 
Moreover, by shifting the phase of $\bfw_n$ appropriately, namely, by replacing $\bfw_n(x,t)$ with $\bfw_n(x,t+\delta_n)$, for suitable $\delta_n$, we may impose
\be
\int_0^{\frac{2\pi}{\omega_n}}\langle\bfw_n,\bfv_{1n}\rangle=\sigma_n\in\real\,,\ \ \int_0^{\frac{2\pi}{\omega_n}}\langle\bfw_n,\bfv_{2n}\rangle=0\,,\ \ \mbox{for all $n\in\nat$.}
\eeq{3.20}


We next observe the following.\\